\documentclass[fleqn,11pt]{article}
\usepackage{amsmath,amssymb,amsthm,xcolor,framed,esint,dsfont}
\usepackage{graphicx}
\usepackage[english]{babel}
\usepackage[margin=3cm]{geometry}
\usepackage{csquotes}
\usepackage{hyperref}

\def\XXint#1#2#3{{\setbox0=\hbox{$#1{#2#3}{\int}$}
		\vcenter{\hbox{$#2#3$}}\kern-.5\wd0}}

\newcommand{\PI}{\mathcal{P}}

\newcommand{\Bs}{\mathsf{B}}

\newcommand{\lt}{\left}
\newcommand{\rt}{\right}

\newcommand{\nn}{\nonumber}

\newcommand{\lm}{\lambda}

\newcommand{\qd}{\quad}

\newcommand{\ep}{\epsilon}

\newcommand{\wt}{\widetilde}

\newcommand{\GI}{\mathcal{G}}

\newcommand{\HI}{\mathcal{H}}
\newcommand{\UI}{\mathcal{U}}

\newcommand{\na}{\nabla}
\newcommand{\ti}{\tilde}
\newcommand{\la}{\langle}
\newcommand{\ra}{\rangle}

\newcommand{\BI}{\mathcal{B}}

\newcommand{\R}{{\mathbb R}}
\newcommand{\Z}{{\mathbb Z}}

\newcommand{\dv}{\nabla\cdot}
\def\Xint#1{\mathchoice
	{\XXint\displaystyle\textstyle{#1}}%
	{\XXint\textstyle\scriptstyle{#1}}%
	{\XXint\scriptstyle\scriptscriptstyle{#1}}%
	{\XXint\scriptscriptstyle\scriptscriptstyle{#1}}%
	\!\int}

\def\XXint#1#2#3{{\setbox0=\hbox{$#1{#2#3}{\int}$}
		\vcenter{\hbox{$#2#3$}}\kern-.5\wd0}}

\newcommand{\e}{\varepsilon}
\DeclareMathOperator{\supp}{supp}
\DeclareMathOperator{\dist}{dist}

\newtheorem{thm}{Theorem}
\newtheorem{prop}[thm]{Proposition}
\newtheorem{lem}[thm]{Lemma}

\newtheorem{rem}[thm]{Remark}

\title{
{On a generalized Aviles-Giga
functional: compactness, zero-energy states, regularity estimates and
energy bounds}
}

\date{}
\author{Xavier Lamy\footnote{Institut de Math\'ematiques de Toulouse, UMR 5219, Universit\'e de Toulouse, CNRS, UPS
IMT, F-31062 Toulouse Cedex 9, France. Email: Xavier.Lamy@math.univ-toulouse.fr}
\and Andrew Lorent\footnote{Department of Mathematical Sciences, University of Cincinnati, Cincinnati, OH 45221, USA. Email: lorentaw@uc.edu} 
\and Guanying Peng\footnote{Department of Mathematical Sciences, Worcester Polytechnic Institute, Worcester, MA 01609, USA. Email: gpeng@wpi.edu}}

\begin{document}

\maketitle

\begin{abstract}
Given any strictly convex norm $\|\cdot\|$ on $\mathbb R^2$ that is $C^1$ in $\R^2\setminus\{0\}$,
we study the generalized Aviles-Giga functional 
\[I_{\epsilon}(m):=\int_{\Omega} \left(\epsilon \left|\nabla m\rt|^2 + \frac{1}{\epsilon}\left(1-\|m\|^2\right)^2\right)    \, dx,\] 
for $\Omega\subset\mathbb R^2$ and $m\colon\Omega\to\mathbb R^2$ satisfying $\nabla\cdot m=0$.
Using, as in the euclidean case $\|\cdot\|=|\cdot|$, the concept of entropies for the limit equation $\|m\|=1$, $\nabla\cdot m=0$, we obtain the following. First, we prove compactness in $L^p$ of sequences of bounded energy. 
Second, we prove rigidity of zero-energy states (limits of sequences of vanishing energy), generalizing and simplifying a result by Bochard and Pegon.
Third, we obtain optimal regularity estimates for limits of sequences of bounded energy, in terms of their entropy productions. 
Fourth, in the case of a limit map in $BV$, we show that lower bound provided by entropy productions and upper bound provided by one-dimensional transition profiles are of the same order.
The first two points are   analogous to what is known in the euclidean case $\|\cdot\|=|\cdot|$, and the last two points are sensitive to the anisotropy of the norm $\|\cdot\|$.
\end{abstract}

\section{Introduction}

The Aviles-Giga functional
\begin{align*}
AG_\e(u)=\int_\Omega \left( \e |\nabla^2 u|^2 +\frac 1\e (1-| \nabla u|^2)^2\right)\, dx,\quad\Omega\subset\R^2,\qd u\colon\Omega\to\R,
\end{align*}
 is a second order functional that (subject to appropriate boundary conditions) models phenomena from thin film blistering to smectic liquid crystals, and is also a natural higher order generalization of the Cahn-Hilliard functional. 
 The  conjecture on  the $\Gamma$-limit of the Aviles-Giga energy, 
 which roughly states that the energy concentrates on a one-dimensional jump set as $\e\to 0$,
  has attracted a great deal of attention, yet remains open; see for example  \cite{avilesgig,avgig1,ADM,mul2,ottodel1, contidel, ark}.
  
 The second term in the  Aviles-Giga functional penalizes values of the divergence-free vector field $m=\nabla^\perp u $ that are far from the euclidean unit circle $\mathbb S^1\subset\R^2$.
In the present work we continue the study, initiated in \cite{boch}, of a generalized Aviles-Giga functional where $\mathbb S^1$ is replaced by the unit circle of  a more general norm on $\R^2$. Specifically, we let $\|\cdot\|$ be a strictly convex norm on $\mathbb R^2$ that is $C^1$ in $\R^2\setminus\{0\}$ (strictly convex $C^1$ norm for simplicity), and consider the generalized Aviles-Giga functional
\begin{align}\label{eq:Ieps}
&I_\e(m)= I_\e(m;\Omega)=\int_\Omega \left(\e |\nabla m|^2 +\frac{1}{\e} (1-\|m\|^2)^2\right)\, dx,\\
& \Omega\subset\R^2,\quad m\colon\Omega\to\R^2, \quad \nabla\cdot m=0\text{ in }\mathcal D'(\Omega).
\nonumber
\end{align}
Here the constraint $\nabla\cdot m=0$ is equivalent to $m=\nabla^\perp u$ if the domain $\Omega$ is simply connected, so $I_\e$ can effectively be seen as a second order functional generalizing the Aviles-Giga functional. In \cite{boch} Bochard and Pegon 
obtain some preliminary results on the characterization of zero-energy states of $I_\e$
(limits of sequences of asymptotically vanishing energy).  In this work we carry out a rather comprehensive analysis of this generalized Aviles-Giga functional. Our goal is to investigate to which extent the results and methods that have been developed for the classical Aviles-Giga functional can be extended to this more general setting. 
In doing so, we hope to shed some light on what parts of the theory are contingent on specific algebraic properties of $\mathbb S^1$, and what parts are more flexible. 
Similar generalized Aviles-Giga functionals have also been studied in \cite{ignat-monteil-20}, 
with a focus on symmetry properties of entire critical points.
Here we concentrate on four aspects: 
\begin{itemize}
\item compactness in $L^p$ and energy lower bounds for sequences of bounded energy; 
\item characterization of zero-energy states; 
\item optimal regularity estimates for limits of sequences of bounded energy;
\item comparison of upper and lower bounds for sequences converging to a map of bounded variation ($BV$). 
\end{itemize}
For the first two aspects  we obtain complete generalizations of the analogous results in the classical case. 
For the last two aspects, our results demonstrate the effects  induced by possible anisotropy and degenerate convexity of $\|\cdot\|$.

A central tool, introduced in \cite{mul2} for the classical Aviles-Giga functional, is the notion of entropies, imported from scalar conservation laws. Formally (and this is   justified by the compactness result), limits of sequences of bounded energy should satisfy the generalized Eikonal equation
\begin{align}\label{eq:geneikon}
\|m\|=1\text{ a.e.},\qquad\nabla\cdot m=0\text{ in }\mathcal D'(\Omega).
\end{align}
Writing locally the unit circle $\partial\Bs=\lbrace z\in\R^2\colon \|z\|=1\rbrace$ as the graph of a convex function $f$, this equation can formally be rewritten as the scalar conservation law 
\begin{align}\label{eq:scl}
\partial_t u +\partial_x f(u) =0.
\end{align}
In direct analogy with the  entropy-entropy flux pairs for this scalar conservation law, entropies for the generalized Eikonal equation \eqref{eq:geneikon} are $C^1$ maps $\Phi\colon\partial \Bs\to\R^2$ with the property that $\nabla\cdot\Phi(m)=0$ for any smooth solution $m$ of \eqref{eq:geneikon}.  For weak solutions, the distributions   $\nabla\cdot\Phi(m)$, called \textit{entropy productions}, encode the presence of singularities and can therefore be used to understand compactness and regularity properties. The key property used in \cite{mul2} is that, in the classical case $\|\cdot\|=|\cdot|$, entropy productions are controlled by the energy. This provides compactness \cite{ADM,mul2}, and an energy lower bound. Thanks to the strict convexity of $\|\cdot\|$, this analysis can be adapted to our generalized setting; see Theorem~\ref{t:comp} and Proposition~\ref{p:lowerbound}.

A further consequence of the energy lower bound is that zero-energy states, that is, limits of sequences $\{m_n\}$ such that $I_{\e_n}(m_n)\to 0$, have vanishing entropy productions. This is exploited in \cite{otto} for the classical Aviles-Giga functional to obtain a kinetic equation which roughly speaking  ensures that zero-energy states are, in a weak way, constant along characteristics. As a consequence, zero-energy states in the classical case $\|\cdot\|=|\cdot|$ are shown in \cite{otto} to be locally Lipschitz outside a locally finite set of singular points, and around each singular point they must coincide with a vortex $m(x)=\pm ix/|x|$. In \cite{boch} this rigidity result is generalized (with appropriate modifications) to $I_\e$ associated with any $C^1$ norm $\|\cdot\|$ of power type $p$ for some $p\in [2,\infty)$ (a quantitative form of strict convexity, see Remark~\ref{r:zeroenergytypep}). Here we extend this further to $I_\e$ associated with any strictly convex $C^1$ norm (see Theorem~\ref{T:zeroenergy}) using an elementary argument that reduces it to the classical case $\|\cdot\|=|\cdot|$.

Finite-energy states, that is, limits of sequences of bounded energy, can have a much more complicated structure. 
The energy lower bound ensures that entropy productions are finite Radon measures, and a central question to solve the $\Gamma$-convergence conjecture for the classical Aviles-Giga functional is whether these measures are concentrated on a one-dimensional rectifiable set. Substantial progress on that question has been made in \cite{ottodel1,marconi21ellipse} but it remains open. 
For scalar conservation laws \eqref{eq:scl} with $f$ uniformly convex (Burgers' equation), this rectifiability property has recently been proved in \cite{marconi20burgers}. 
The results of \cite{ottodel1} and \cite[Proposition~1.7]{marconi21ellipse} can likely be generalized to the class of energy functionals \eqref{eq:Ieps} associated with any strictly convex $C^1$ norm $\|\cdot\|$ (using the kinetic formulation obtained in Lemma~\ref{l:kin}), 
but here we don't address that question and concentrate instead on optimal regularity estimates for solutions of the generalized Eikonal equation \eqref{eq:geneikon} whose entropy productions are locally finite Radon measures. 
In the classical case $\|\cdot\|=|\cdot|$, it is proved in \cite{GL} (adapting an argument of \cite{golseperthame13} for scalar conservation laws) that such solutions must locally have the Besov regularity $B^{\frac 13}_{3,\infty}$, i.e. $\sup_h |h|^{-\frac 13}\|m-m^h\|_{L^3_{loc}}<\infty$ where $m^h=m(\cdot+h)$.
Moreover this estimate is \emph{strongly optimal} in the sense that it is \emph{equivalent} to  entropy productions being locally finite Radon measures. 
In the general case, the coercivity provided by the strict convexity of the norm $\|\cdot\|$ depends on the direction $z$ on its unit circle $\partial\Bs$, and optimal estimates must take that into account. We prove therefore a regularity estimate of the form 
$\sup_h |h|^{-1}\| \Pi(m,m^h)\|_{L^1_{loc}}<\infty$ for some function $\Pi\colon\partial\Bs\times\partial\Bs\to [0,\infty)$ that is sensitive to the anisotropy of $\|\cdot\|$, 
and show that it is strongly optimal (equivalent to entropy productions being locally finite Radon measures) at least when the norm $\|\cdot\|$ is analytic; see Theorems~\ref{thmbes} and \ref{thmbesrev}. 
(For a norm $\|\cdot\|$ of power type $p$ convexity  this estimate implies in particular Besov regularity agreeing with the results of \cite{golseperthame13} for scalar conservation laws \eqref{eq:scl} when the flux $f$ has degenerate convexity; see Remark~\ref{R6}.)
Furthermore, if $\|\cdot\|$ is merely $C^1$ then the quantity
 $\sup_h |h|^{-1}\| \Pi(m,m^h)\|_{L^1_{loc}}$  is comparable to the total entropy production when $m$ is $BV$, hinting that the regularity estimate could be strongly optimal for all strictly convex $C^1$ norms $\|\cdot\|$.

The $\Gamma$-convergence of the classical Aviles-Giga functional in the $BV$ setting is well understood \cite{ADM, contidel, ark}. For a solution $m$ of the generalized Eikonal equation \eqref{eq:geneikon} which is $BV$, an upper bound can be obtained for the minimal energy of approximating sequences $m_n\to m$ by pasting optimal one-dimensional transitions along the jump set $J_m$ \cite{contidel,ark,ark2}. 
In the classical case $\|\cdot\|=|\cdot|$, this upper bound happens to coincide with the lower bound provided by a particular class of entropy productions \cite{ADM}, thus characterizing the $\Gamma$-limit at $BV$ maps $m$. 
This perfect agreement of entropy lower bound and 1D upper bound is very likely linked to specific algebraic properties of the euclidean norm $|\cdot|$ (as are the symmetry results of \cite{ignat-monteil-20}). 
In fact it is known \cite[\S~4]{kohn} that for general $\|\cdot\|$ optimal transition profiles may not be one-dimensional, and in that case the 1D upper bound is strictly larger than any lower bound (see \cite{poliakovsky13,poliakovsky15} for more results related to such issues).
It is however interesting to find out whether these two bounds (the entropy lower bound and the 1D upper bound) are of the same order of magnitude, or can instead be very far apart. 
Like optimal regularity estimates, this question is sensitive to the possibly  anisotropic behavior of $\|\cdot\|$. We prove that these upper and lower bounds do agree up to a multiplicative constant; see Theorem~\ref{t:lowerupper}.

In the rest of this introduction we present the precise statements of our results. 
In Section~\ref{s:entropyprod} we derive some useful properties of the entropies in our generalized setting. In Section~\ref{s:compactness} we prove the compactness result. In Section~\ref{s:zeroenergy} we prove the rigidity of zero-energy states. In Section~\ref{s:reg} we prove regularity estimates  for finite-energy states and their optimality. And in Section~\ref{s:boundsBV} we compare upper and lower bounds for $BV$ limits.

\subsection{Notations and assumptions}

Let $\Omega\subset\R^2$ be a bounded open set and $\|\cdot\|$ be a strictly convex $C^1$ norm on $\R^2$ unless otherwise specified. We denote by $\Bs=\left\lbrace z\in\R^2\colon \|z\|<1\right\rbrace$ the open unit disk for the norm $\|\cdot\|$. The properties of $\|\cdot\|$ are equivalent to strict convexity of $\Bs$ and $\partial \Bs$ being a $C^1$ manifold.
Without loss of generality, we assume that $\partial\Bs$ has length $2\pi$, and let
 $\gamma:\R/2\pi\Z\rightarrow \partial \Bs$ be the counterclockwise arc-length parametrization of $\partial\Bs$ (unique up to translation of the variable). By assumption, $\gamma\in C^1(\R/2\pi\mathbb{Z};\R^2)$. In many places we identify $\mathbb{R}^2$ with $\mathbb{C}$ and in particular we let $i$ denote the counterclockwise rotation by $\frac{\pi}{2}$. 
  We will  use the symbols $\lesssim$ and $\gtrsim$ to denote inequality up to a multiplicative constant that depends only on $\Bs$.

\subsection{Compactness and lower bound} 

Our first result generalizes the compactness result obtained independently in  \cite[Theorem~3.3]{ADM} and  \cite[Proposition~1]{mul2} for the Aviles-Giga functional.
\begin{thm} 
	\label{t:comp} 
	 Suppose  the sequence $\{m_{n}\}\subset W^{1,2}(\Omega;\R^2)$  satisfies $\nabla\cdot m_n=0$ and
	\begin{equation*}
		\sup_{n} I_{\ep_n}(m_{n})<\infty.
	\end{equation*}
	Then $\{m_{n}\}$ is precompact in $L^2(\Omega)$.
\end{thm}

As explained above, this compactness result relies heavily on the notion of entropies for the generalized Eikonal equation  
\begin{align}\tag{\ref{eq:geneikon}}
\|m\|=1\text{ a.e.,}\quad \nabla\cdot m=0\text{ in }\mathcal D'(\Omega).
\end{align}
Equivalently, the first constraint $\|m\|=1$ means that $m$ takes values into $\partial\Bs$.
Entropies for this equation are $C^1$ maps $\Phi\colon\partial \Bs\to\R^2$ such that, if $m$ is a $C^1$ solution of  \eqref{eq:geneikon}, then $\Phi(m)$ is also divergence-free $\nabla\cdot\Phi(m)=0$. It is a lengthy but straightforward exercise to see that this is equivalent to requiring that, for all $\theta\in\R$,
\begin{align*}
\frac{d}{d\theta}\Phi(\gamma(\theta)) \text{ is tangent to }\partial\Bs\text{ at }\gamma(\theta).
\end{align*}
For a weak solution $m$ of \eqref{eq:geneikon}, the entropy production $\nabla\cdot\Phi(m)$ is in general not zero, and encodes the presence of singularities. The proof of 
Theorem~\ref{t:comp} relies on the control of entropy productions provided by the  energy. This control is possible for  regular enough  entropies: we define
\begin{align}\label{eq:ENT}
\mathrm{ENT}=\Big\lbrace \Phi\in C^1(\partial\Bs;\R^2)\colon &\frac{d}{d\theta}\Phi(\gamma(\theta))=\lambda_\Phi(\theta)\gamma'(\theta)
\nonumber\\
&\text{for some function }\lambda_\Phi\in C^1(\R/2\pi\Z)\Big\rbrace.
\end{align}
The control of entropy productions used to establish compactness also provides a lower bound for the energy. From this point on all entropies for equation \eqref{eq:geneikon}  in statements and proofs will be taken to be the ones from $\mathrm{ENT}$.
\begin{prop}\label{p:lowerbound}
Let $m\colon\Omega\to\R^2$ be such that $m=\lim_{n\to\infty} m_n$ in $L^2(\Omega)$ for some sequence $\{m_{n}\}\subset W^{1,2}(\Omega;\R^2)$ with $\nabla\cdot m_n=0$ and $\sup_n I_{\e_n}(m_n) <\infty$. Then $m$ satisfies the generalized Eikonal equation \eqref{eq:geneikon}, its entropy productions satisfy $\nabla\cdot\Phi(m)\in\mathcal{M}(\Omega)$ for all $\Phi\in \mathrm{ENT}$, and they provide the lower bound
\begin{align}\label{eq:lowerbound}
\left(\bigvee_{\|\lambda_\Phi'\|_{\infty}\leq 1}|\nabla\cdot\Phi(m)|\right)(U) \leq C_0 \,\liminf_{n\to\infty} I_{\e_n}(m_n;U),
\end{align}
for any open subset $U\subset\Omega$ and some constant $C_0>0$ depending only on $\Bs$. Here $\bigvee$ denotes the lowest upper bound measure \cite[Definition~1.68]{ambrosio} of a family of measures.
\end{prop}

\begin{rem} The hypothesis that $\|\cdot \|$ is strictly convex is necessary for Theorem~\ref{t:comp}: Suppose that $\partial \Bs$ contains a line segment $\lt[\zeta_0,\zeta_1\rt]$ then without loss of generality we can assume $\zeta_0=e_1+\delta e_2$ and $\zeta_1=e_1-\delta e_2$. 
Setting 
$m_\e(x)=e_1 + \delta \sin \lt(x_1/\sqrt{\e}\rt) e_2$, then $\na \cdot m_\e=0$ and $\|m_\e\|=1$ everywhere in $\Omega$. Thus $\sup_{\e>0} I_{\e}\lt(m_\e\rt)<\infty$, but $m_\e$ converges weakly to $\ti{m}\equiv e_1$ in $L^p$ as $\e\to 0$ and $\|m_\e-\ti{m}\|_{L^p}\gtrsim \delta$ for all $\e>0$ and all $p\geq 1$.
\end{rem}

\subsection{Zero-energy states}

As stated previously, Jabin, Otto and Perthame showed in  \cite[Theorem~1.1]{otto} that zero-energy states of the Aviles-Giga functional are  rigid. 
This result
 has several interesting implications \cite{ottodel1,DeI,LP,llp}. 
It is proved in two steps: first, zero-energy states have vanishing entropy productions and satisfy as a consequence the kinetic equation $e^{it}\cdot\nabla_x\mathbf 1_{m(x)\cdot e^{it}>0}=0$, which expresses in a weak way the fact that $m$ is constant along characteristics of the classical Eikonal equation; second, solutions of this kinetic equation are shown to be rigid.
In \cite{boch}, Bochard and Pegon generalize the second step to solutions of the kinetic equation $\gamma'(t)\cdot\mathbf 1_{m(x)\cdot i\gamma(t)>0}=0$ naturally associated with the generalized Eikonal equation \eqref{eq:geneikon}, under the assumption that the $C^1$ norm $\|\cdot\|$ is of power type $p$ (see Remark~\ref{r:zeroenergytypep}). They do not however prove the first step, namely that zero-energy states of $I_\e$ satisfy this kinetic equation.
Here we do establish that missing step, and generalize their rigidity result to any strictly convex $C^1$ norm $\|\cdot\|$, with a somewhat more direct proof.

\begin{thm}
	\label{T:zeroenergy}
	Let $m:\Omega\to \R^2$ be such that $m=\lim_{n\to\infty} m_n$ in $L^2(\Omega)$, where the sequence $\{m_n\}\subset W^{1,2}(\Omega;\R^2)$ satisfies $\na\cdot m_n=0$ and
	\begin{equation*}
		\lim_{n\to\infty}I_{\ep_n}(m_n) = 0.
	\end{equation*}
	Then $m$ is continuous outside a locally finite set of singular points. For every singular point $x_0$, there exists $\beta\in\{-1,1\}$ such that in any convex neighborhood $\UI$ of $x_0$, we have $m(x)=\beta V_{\Bs}(i(x-x_0))$, where $V_{\Bs}(\cdot):=\na \|\cdot \|_{*}$ is the vortex associated to $\|\cdot\|$ 
	and $\|\cdot \|_{*}$ is the dual norm of $\|\cdot \|$.
\end{thm}

\begin{rem}\label{r:zeroenergytypep}
Our proof also recovers the result, obtained in \cite{boch}, that if the $C^1$
norm $\|\cdot\|$  is of power type $p$  for some $p\geq 2$, that is,
\begin{align*}
1-\left\|\frac{x+y}{2}\right\| \geq K \|x-y\|^p\qquad\forall x,y\in\partial \Bs,
\end{align*}
for some constant $K>0$,
 then $m$ is locally $\frac{1}{p-1}$-H\"{o}lder outside a locally finite set of singular points (see the end of Section \ref{s:zeroenergy}). 
\end{rem}

\subsection{Optimal regularity estimates}

Proposition \ref{p:lowerbound} motivates the study of \emph{finite-entropy solutions} of the generalized Eikonal equation, i.e. solutions of \eqref{eq:geneikon} satisfying $\na\cdot\Phi(m)\in\mathcal{M}_{loc}(\Omega)$ for all $\Phi\in\mathrm{ENT}$.
 We present here regularity estimates for these solutions, that are strongly optimal in the sense that a converse estimate is valid: regularity implies locally finite entropy productions. 
 In the context of scalar conservation laws, this type of optimality is related to ``Onsager conjecture-type'' statements: see e.g.  \cite{bardos19} where the authors investigate minimal regularity requirements that are sufficient to ensure that entropy productions vanish.

In the classical case $\|\cdot\|=|\cdot|$ it was shown in \cite{GL} that finite-entropy solutions coincide  with solutions of \eqref{eq:geneikon} that live in the Besov space $B^{1/3}_{3,\infty,loc}(\Omega)$: such Besov estimates are strongly optimal. 
This was obtained by adapting methods of \cite{golseperthame13} for scalar conservation laws \eqref{eq:scl} with convex flux $f$. 
The $1/3$ order of regularity is valid for uniformly convex fluxes $f$ and
an example in \cite{delelliswestdickenberg03} had also demonstrated its optimality, in a different sense than the one we wish to study here: there exist finite-entropy solutions which don't have a better order of regularity.

For fluxes with degenerate convexity, quantified by the inequality
\begin{align*}
f'(v)-f'(w)\gtrsim \lt|v-w\rt|^{p-1}\qquad \forall v>w,
\end{align*}
for some $p\geq 2$, the regularity obtained in \cite{golseperthame13} is $B^{\frac 1{p+1}}_{p+1,\infty,loc}$. This applies for instance to $f(w)=|w|^p$, and is shown to be  optimal in \cite[Proposition~3.2]{delelliswestdickenberg03}, again  in the sense that there exist finite-entropy solutions which don't have a better order of regularity.
However it is clear (considering solutions whose values stay away from the point $w=0$ at which convexity degenerates) that this Besov regularity does not provide a converse estimate: it is not strongly optimal.

Here, following \cite{GL} we adapt the methods of \cite{golseperthame13} to the generalized Eikonal equation \eqref{eq:geneikon} in order to obtain regularity estimates that take into account the anisotropy of $\|\cdot\|$, and in particular the fact that the convexity of $\|\cdot\|$ may degenerate differently in different directions. For a precise statement, we introduce the (unique up to an additive constant integer multiple of $2\pi$) continuous function   $\alpha\colon\R\to\R$  such that
\begin{align*}
\gamma'(\theta)=e^{i\alpha(\theta)}\qquad\forall \theta\in\R.
\end{align*}
The strict convexity of $\Bs$ ensures that this function $\alpha$ is increasing, and the symmetry of $\Bs$ implies $\alpha(t+\pi)=\alpha(t)+\pi$ for all $t\in\R$. We define a function $\Pi\colon\partial \Bs\times \partial \Bs\to [0,\infty)$ by  
\begin{align}
\label{eq:Lambda}
\Pi\lt(\gamma(\theta_1),\gamma(\theta_2)\rt) =\int_{\theta_1}^{\theta_2} \int_{\theta_1}^{\theta_2} |\alpha(t)-\alpha(s)|\, dt ds \qquad\text{ for }|\theta_1-\theta_2|\leq \pi.
\end{align}
Using this function $\Pi$ as a ``metric'' for the increments, we have the following regularity estimate for finite-entropy solutions of \eqref{eq:geneikon}.

\begin{thm}
\label{thmbes} Let $m$ satisfy  the generalized Eikonal equation  \eqref{eq:geneikon}. Suppose
\begin{align}
\label{eqbes1}
\na \cdot \Phi(m)\in \mathcal{M}_{loc}(\Omega)\qd\text{ for all }\Phi\in \mathrm{ENT},
\end{align}
then 
\begin{align}
\label{eqbes2}
\sup_{\lt|h\rt|<\mathrm{dist}(\Omega',\partial \Omega)} \frac{1}{\lt|h\rt|}\int_{\Omega'} \Pi\lt(m(x+h),m(x)\rt)\, dx<\infty\qd \text{ for any }\Omega'\subset \subset \Omega.
\end{align}
%
\end{thm}

\begin{rem}
\label{R6}
If $\|\cdot\|$ is of power type $p$ for some $p\geq 2$ (see Remark~\ref{r:zeroenergytypep}), the estimate \eqref{eqbes2} directly implies that $m\in B^{\frac{1}{p+1}}_{p+1, \infty, loc}(\Omega)$, as explained in Remark~\ref{remppower}. This corollary is analogous to the regularity results obtained in \cite[Theorem~4.1]{golseperthame13} for convex scalar conservation laws.
\end{rem}

The main interest of this regularity estimate is that the ``metric'' $\Pi$ is sensitive enough to the local convexity of $\partial\Bs$ to ensure the validity of a converse estimate, at least when $\|\cdot\|$ is analytic  in $\R^2\setminus \lbrace 0\rbrace$ (or equivalently $\partial\Bs$ is analytic):

\begin{thm}
\label{thmbesrev}
Let $m$ satisfy \eqref{eq:geneikon}. Assume that the strictly convex norm $\|\cdot\|$ is  analytic  in $\R^2\setminus\lbrace 0\rbrace$, then \eqref{eqbes2} implies \eqref{eqbes1}.
\end{thm}

\begin{rem}
Note that Theorem~\ref{thmbesrev} applies in particular to $\|\cdot\|=\|\cdot\|_{\ell^p}$ for any $1<p<\infty$ (see Remark \ref{rem:lp-sphere}). 
\end{rem}

We don't know whether the analyticity assumption on $\|\cdot\|$ is necessary for the validity of the converse estimate  \eqref{eqbes2} implying \eqref{eqbes1}. An indication that it might not be needed is given by the following.

\begin{thm}\label{t:regconvBV}
Let $\|\cdot\|$ be a strictly convex $C^1$ norm on $\R^2$ and let $m\in BV(\Omega; \R^2)$ satisfy \eqref{eq:geneikon}. Then for any open subset $\Omega'\subset\subset\Omega$ we have
\begin{align*}
\left(\bigvee_{\|\lambda_\Phi'\|_\infty\leq 1} |\nabla\cdot\Phi(m)|\right)(\Omega')
\leq C_0 \,
\sup_{|h|<\dist(\Omega',\partial\Omega)}\frac{1}{|h|}\int_{\Omega'}\Pi(m(x+h),m(x))\, dx,
\end{align*}
for some absolute constant $C_0>0$.
\end{thm}

Note that for a $BV$ solution of \eqref{eq:geneikon} as in Theorem~\ref{t:regconvBV}, both the entropy productions and the quantity appearing in the regularity estimate \eqref{eqbes2} are finite. Here the point is that the latter controls the former, without any further regularity assumption on the norm $\|\cdot\|$.

%

%
%

\subsection{Comparison of upper and lower bounds}

For general maps $m$, finding an upper bound that matches (at least up to a multiplicative constant) the lower bound of Proposition~\ref{p:lowerbound} is a famously hard problem even in the classical case. However, when the limiting solution $m$ of \eqref{eq:geneikon} additionally belongs to $BV(\Omega;\R^2)$, then it is known \cite{ark2} that  an upper bound (in the sense of $\Gamma$-convergence) is obtained by pasting optimal one-dimensional transitions at scale $\e$ along the jump set $J_m$. Specifically, for any solution $m\in BV(\Omega;\R^2)$ of the generalized Eikonal equation \eqref{eq:geneikon} and any smooth simply connected open subset $U\subset\Omega$, there exists a sequence $m_\e\to m$ in $L^p(U;\R^2)$ for $1\leq p <\infty$, such that
\begin{align}\label{eq:upperboundBV}
\limsup_{\e\to 0}I_{\e}(m_\e;U) \leq \int_{U\cap J_m} \mathrm{c}^{\mathrm{1D}}(m^+,m^-)\, d\mathcal H^1,
\end{align}
where $m^\pm$ are the traces of $m$ along $J_m$, and $\mathrm{c}^{\mathrm{1D}}\colon\partial\Bs\times\partial\Bs\to [0,\infty)$ is given by
\begin{align}\label{eq:c1D}
\mathrm{c}^{\mathrm{1D}}(z^+,z^-)
&
=2\left|\int_{z^{-}\cdot i\nu}^{z^{+}\cdot i\nu} \lt(1-\|a \nu + s i \nu \|^2\rt) ds\right|,\\
\nu&=i\frac{z^+-z^-}{|z^+-z^-|},\qquad a=z^+\cdot\nu =z^-\cdot \nu.\nonumber
\end{align}
Here the unit vector $\nu$ represents a normal vector to the jump set $J_m$ at a jump between $z^+$ and $z^-$. The divergence-free constraint $\nabla\cdot m=0$ forces $\nu$ to satisfy $(z^+-z^-)\cdot\nu=0$, and this characterizes $\nu$ up to a sign. Note that it is known that the upper bound provided by one-dimensional profiles will in general not be optimal \cite[\S~4]{kohn} (see \cite{poliakovsky13,poliakovsky15} for a discussion of optimal upper bounds), but here we are only interested in optimality up to a multiplicative constant.

We wish to compare this 1D upper bound to the lower bound provided by the entropy productions in Proposition~\ref{p:lowerbound}. For
 a solution $m$ of \eqref{eq:geneikon}  which additionally belongs to $BV(\Omega;\R^2)$, the $BV$ chain rule implies that the entropy productions are absolutely continuous with respect to $\mathcal H^1_{\lfloor J_m}$. Thanks to \cite[Remark~1.69]{ambrosio}, the resulting lowest upper bound measure is also absolutely continuous  with respect to $\mathcal H^1_{\lfloor J_m}$, and (see Lemma \ref{l:cENT}) we have
\begin{align*}
\bigvee_{\|\lambda'_\Phi\|_{\infty}\leq 1}|\nabla\cdot\Phi(m)| & =\mathrm{c}^{\mathrm{ENT}}(m^+,m^-)\,\mathcal H^1_{\lfloor J_m},
\end{align*}
where $m^\pm$ are the traces of $m$ along $J_m$, and the jump cost $\mathrm{c}^{\mathrm{ENT}}\colon \partial\Bs\times\partial\Bs\to\ [0,\infty)$ is given by
\begin{align*}
\mathrm{c}^{\mathrm{ENT}}(z^+,z^-)=\sup_{\|\lambda'_\Phi\|_{\infty}\leq 1}
\int_{\theta^-}^{\theta^+}\lambda_\Phi(s)\,\gamma'(s)\cdot\nu\, ds\qquad\text{for }z^\pm=\gamma(\theta^\pm),
\end{align*}
and $\nu$ as in \eqref{eq:c1D}.
The value of the last integral does not depend on the choices of $\theta^\pm$ modulo $2\pi$, because  the definition of $\lambda_\Phi$ in \eqref{eq:ENT} implies  $\int_{\R/2\pi\Z}\lambda_\Phi(s)\gamma'(s)\, ds=0$ for any entropy $\Phi\in \mathrm{ENT}$.
In other words, for a $BV$ map $m$, the lower bound \eqref{eq:lowerbound} becomes
\begin{align}\label{eq:lowerboundBV}
\liminf_{n\to\infty} I_{\e_n}(m_n;U)
\gtrsim
 \int_{J_m\cap U} \mathrm{c}^{\mathrm{ENT}}(m^+,m^-)\, d\mathcal H^1.
\end{align}
We show that these lower and upper bounds  \eqref{eq:lowerboundBV} and \eqref{eq:upperboundBV}  for $BV$ maps $m$ are comparable:
\begin{thm}\label{t:lowerupper}
There exists a constant $C_0>0$ depending only on $\Bs$ such that
\begin{align*}
C_0^{-1}\,\mathrm{c}^{\mathrm{1D}}(z^+,z^-)\leq \mathrm{c}^{\mathrm{ENT}}(z^+,z^-)
\leq C_0\, \mathrm{c}^{\mathrm{1D}}(z^+,z^-) 
\end{align*}
for all $z^\pm\in\partial\Bs$.
\end{thm}

\section{Entropy productions}\label{s:entropyprod}

In this section we compute entropy productions of divergence-free $m\in W^{1,2}$, and as a direct consequence we prove Proposition \ref{p:lowerbound}. Since $\Bs$ is convex and centered, for any $z\in \R^2\setminus\lbrace 0\rbrace$ there is a unique $(r,\theta)\in (0,\infty)\times\R/2\pi\Z$ such that $z=r\gamma(\theta)$. In order to make use of classical polar coordinates, we introduce the bijection $X\colon\R^2\to\R^2$ given by
\begin{align}\label{eq:X}
X(r e^{i\theta})= r \gamma(\theta)\qquad\forall r\geq 0,\;\theta\in\R.
\end{align}
The map $X$ is $C^1$ in $\R^2\setminus\lbrace 0\rbrace$, and its jacobian determinant is
\begin{align*}
\det (\nabla X(re^{i\theta})) =i\gamma(\theta)\cdot\gamma'(\theta) \geq\alpha_0>0,
\end{align*}
where $\alpha_0$ is the radius of the largest euclidean ball contained in $\overline\Bs$. This last inequality follows from the convexity of $\Bs$: for any $z\in \overline\Bs$ we have $(z-\gamma(\theta))\cdot i\gamma'(\theta) \geq 0$, and applying this to $z=-i\alpha_0\gamma'(\theta)$ gives $i\gamma(\theta)\cdot\gamma'(\theta) \geq\alpha_0$.  
As a consequence, $X^{-1}$ is $C^1$ in $\R^{2}\setminus\lbrace 0\rbrace$. Moreover $X$ is a bi-Lipschitz homeomorphism.

In the following, we take $\eta(r)\in C^1([0,\infty))$ so that $0\leq\eta\leq 1$, $\eta\equiv 0$ in $[0,\frac 12]\cup[2,\infty)$ and $\eta(1)=1$. For $\Phi\in \mathrm{ENT}$, define $\widehat\Phi\in C^1(\R^2;\R^2)$ by
\begin{equation}\label{eq:extPhi}
	\widehat\Phi\lt(r\gamma(\theta)\rt) = \eta(r)\Phi\lt(\gamma(\theta)\rt).
\end{equation}

\begin{lem}
	\label{l:extPhi}
	Let $m\in W^{1,2}(\Omega;\R^2)$ satisfy $\dv m=0$. Then for any $\Phi\in \mathrm{ENT}$,
	we have 
	\begin{align} \label{eq:dvextPhim}
		\dv\widehat\Phi(m)=
		\frac 12\Psi(m)\cdot\na\lt(1-\|m\|^2\rt),
	\end{align}
	where
	\begin{equation*}
		\Psi\lt(r\gamma(\theta)\rt)=\frac{\eta\lt(r\rt)\lm_{\Phi}\lt(\theta\rt)}{r^2}\gamma\lt(\theta\rt) - \frac{\eta'\lt(r\rt)}{r}\Phi\lt(\gamma(\theta)\rt)
	\end{equation*}
	and 
	\begin{equation*}
		\lm_{\Phi}(\theta)=\frac{d}{d\theta}\lt(\Phi\lt(\gamma(\theta)\rt)\rt)\cdot \gamma'(\theta).
	\end{equation*}
\end{lem}
\begin{proof}
It suffices to prove \eqref{eq:dvextPhim} for a smooth map $m\colon\Omega\to\R^2$, because we can then approximate a $W^{1,2}$ map $m$ with smooth maps $m_n\to m$ in $W^{1,2}$ and a.e., satisfying in addition $\nabla\cdot m_n=0$, so that \eqref{eq:dvextPhim} passes to the limit in $\mathcal D'(\Omega)$.

It is convenient to change variable in order to use classical polar coordinates: we set
$\widetilde\Phi =\widehat\Phi\circ X$, $\tilde m = X^{-1}(m)$,
so that 
\begin{align}
\label{eqfcc1}
\widetilde\Phi(r e^{i\theta})=\eta(r)\Phi(\gamma(\theta)),
\qquad\widetilde\Phi(\tilde m)=\widehat\Phi(m), \qquad\text{and }\na\cdot X(\ti m)=0. 
\end{align}
And note that $\|X(v)\|=\lt|v\rt|$ for all $v\in \mathbb{R}^2$, and so
\begin{align}
\label{eqfcc2}
 \lt|\ti{m}\rt|=\|m\|.
\end{align}
In the following we perform calculations in the open set $\{m\ne 0\}$. Using that $\Phi \in\mathrm{ENT}$, we have
\begin{align*}
\frac{d}{d\theta}\Phi(\gamma(\theta))=\lambda_{\Phi}(\theta)\gamma'(\theta),
\end{align*}
so computing $D\widetilde\Phi$ using polar coordinates we find
\begin{align*}
D\widetilde\Phi(re^{i\theta})&=\eta'(r)\Phi(\gamma(\theta))\otimes e^{i\theta}
+\frac{\eta(r)}{r}\lambda_{\Phi}(\theta) \,\gamma'(\theta)\otimes ie^{i\theta}.
\end{align*}
Note that for any $v,w\in \mathbb{R}^2$ we have the identity $\mathrm{tr}\lt(v\otimes w Dm\rt)=\lt(\lt(v\cdot \na\rt)m \rt)\cdot w$. So, writing $\tilde m=re^{i\theta}$, we obtain
\begin{align*}
\nabla\cdot\widetilde\Phi(\tilde m)&=\mathrm{tr}(D\widetilde \Phi(\tilde m)D\tilde m)\\
&=\eta'(|\tilde m|)\lt[\lt(\Phi(\gamma(\theta))\cdot\nabla\rt)\tilde m\rt] \cdot \frac{\tilde m}{|\tilde m|}
+\frac{\eta(|\tilde m|)}{|\tilde m|}\lambda_{\Phi}(\theta) \lt[\lt(\gamma'(\theta)\cdot\nabla\rt)\tilde m\rt] \cdot i
\frac{\tilde m}{|\tilde m|}.
\end{align*}
Applying this to $\eta(r)=r$ (the above calculations only require $\eta$ to be $C^1$) and $\Phi(z)=z$ gives in particular
\begin{align*}
\nabla\cdot X(\tilde m)=
\lt[\lt(\gamma(\theta)\cdot\nabla\rt)\tilde m\rt] \cdot \frac{\tilde m}{|\tilde m|}
+ \lt[\lt(\gamma'(\theta)\cdot\nabla\rt)\tilde m\rt] \cdot i
\frac{\tilde m}{|\tilde m|},
\end{align*}
so the previous expression for $\nabla\cdot\widetilde\Phi(\tilde m)$ can be rewritten as
\begin{align*}
\nabla\cdot\widetilde\Phi(\tilde m)&=
\frac{\eta(|\tilde m|)}{|\tilde m|}\lambda_{\Phi}(\theta) \,\nabla\cdot X(\tilde m) \\
&\quad
+
\lt[\lt(\left(\eta'(|\tilde m|)\Phi(\gamma(\theta))
-\frac{\eta(|\tilde m|)}{|\tilde m|}\lambda_{\Phi}(\theta) \gamma(\theta)\right)\cdot\nabla\rt)\tilde m \rt]
\cdot \frac{\tilde m}{|\tilde m|}.
\end{align*}
Using $\nabla\cdot X(\tilde m)=0$ and $\partial_j\ti m\cdot\ti m = \partial_j |\ti m|^2/2$, this becomes
\begin{align*}
\nabla\cdot\widetilde\Phi(\tilde m)&=
\frac 12 \widetilde\Psi(\tilde m)\cdot\nabla (1- |\tilde m|^2),\\
\widetilde\Psi(\tilde m)&=
\frac{\eta(|\tilde m|)}{|\tilde m|^2}\lambda_{\Phi}(\theta) \gamma(\theta)
-
\frac{\eta'(|\tilde m|)}{|\tilde m|}\Phi(\gamma(\theta)).
\end{align*}
The above calculations are valid in $\lbrace m\neq 0\rbrace$, but since $\eta(r)=\eta'(r)=0$ for $0\leq r <1/2$, this last expression makes sense everywhere.
Recalling  from \eqref{eqfcc1}-\eqref{eqfcc2} that $\widetilde\Phi(\tilde m)=\widehat\Phi(m)$, and $|\tilde m|=\|m\|$, setting $\Psi=\wt{\Psi}\circ X^{-1}$ this is exactly the claimed expression \eqref{eq:dvextPhim} for $\nabla\cdot\widehat\Phi(m)$.
\end{proof}

Proposition~\ref{p:lowerbound} is a rather direct consequence of the identity obtained in Lemma~\ref{l:extPhi}.

\begin{proof}[Proof of Proposition~\ref{p:lowerbound}]
Let $m\colon\Omega\to \R^2$ be such that $m=\lim_{n\to\infty} m_n$ in $L^2(\Omega)$ 
for some $\{m_n\}\subset W^{1,2}(\Omega;\R^2)$
 with $\nabla\cdot m_n=0$ and  $\sup_n I_{\e_n}(m_n)<\infty$. 
The fact that $\na\cdot m=0$ in $\mathcal{D}'(\Omega)$ follows from $\nabla\cdot m_n=0$ and $L^2$ convergence. 
The assumption $\sup_n I_{\e_n}(m_n)<\infty$ implies that $\|m_n\|\to 1$ in $L^2(\Omega)$. 
This together with $m_n\to m$ in $L^2(\Omega)$ gives $\|m\|=1$ a.e., and thus $m$ satisfies the generalized Eikonal equation \eqref{eq:geneikon}.

Let $\Phi\in \mathrm{ENT}$ and its extension $\widehat\Phi$ defined in \eqref{eq:extPhi}.
First note that, in order to estimate $\nabla\cdot\Phi(m)$, we may assume without loss of generality that
\begin{align}\label{eq:Phizeroavg}
\int_{\R/2\pi\Z}\lambda_\Phi =0\quad\text{and}\quad \int_{\partial\Bs}\Phi =0.
\end{align}
This is due to the fact that, for any $a\in\R$ and $b\in\R^2$ the entropy given by $\Phi^{a,b}(z)=\Phi(z)+az+b$ for $z\in\partial\Bs$ satisfies $\nabla\cdot\Phi^{a,b}(m)=\nabla\cdot\Phi(m)$ since $\nabla\cdot m=0$, and $\lambda_{\Phi^{a,b}}=\lambda_\Phi + a$. Hence we may choose $a$ such that $\lambda_{\Phi^{a,b}}$ has zero average for any $b\in\R^2$, and $b$ such that $\Phi^{a,b}$ has zero average.

Thanks to Lemma~\ref{l:extPhi}, for any test function $\zeta\in C_c^\infty(\Omega)$ with support inside an open subset $V\subset\Omega$, we have
\begin{align*}
\langle \nabla\cdot\widehat\Phi(m_n),\zeta\rangle 
& = - \frac 12\int_\Omega \zeta\, (1-\|m_n\|^2) \nabla\cdot\Psi(m_n) \, dx\nn\\
& \quad 
-\frac 12 \int_\Omega \Psi(m_n)\cdot\nabla\zeta\, (1-\|m_n\|^2)\, dx \nn\\
& \lesssim \|\nabla\Psi\|_\infty \|\zeta\|_\infty I_{\e_n}(m_n;V) + 
\|\Psi\|_\infty \|\nabla\zeta\|_\infty |V|^{\frac 12} \e_n^{\frac 12} I_{\e_n}(m_n)^{\frac 12},
\end{align*}
so taking the limit $n\to\infty$ we deduce
\begin{align}\label{eq:dvPhimzeta}
\langle\nabla\cdot\Phi(m),\zeta\rangle  \lesssim \|\nabla\Psi\|_\infty \|\zeta\|_\infty \liminf_{n\to\infty}I_{\e_n}(m_n;V).
\end{align}
This implies in particular that $\nabla\cdot\Phi(m)$ is a  finite Radon measure. From the proof of Lemma \ref{l:extPhi}, we have $\Psi=\widetilde\Psi\circ X^{-1}$ with
\begin{align*}
\widetilde\Psi(re^{i\theta})=\frac{\eta\lt(r\rt)\lm_{\Phi}\lt(\theta\rt)}{r^2}\gamma\lt(\theta\rt) - \frac{\eta'\lt(r\rt)}{r}\Phi\lt(\gamma(\theta)\rt).
\end{align*}
Recalling that $X^{-1}$ is Lipschitz we deduce
\begin{align*}
\|\nabla\Psi\|_\infty \lesssim \|\Phi\|_{C^1} + \|\lambda_{\Phi}\|_{C^1}.
\end{align*}
Recall that $\lambda_\Phi$ and $\Phi$ have zero average thanks to \eqref{eq:Phizeroavg} and thus $\|\lambda_{\Phi}\|_{C^1}$ is controlled by  $\|\lambda_{\Phi}'\|_{\infty}$. 
Further, as $(d/d\theta)\Phi(\gamma(\theta))=\lambda_\Phi(\theta)\gamma'(\theta)$, we also have that $\|\Phi\|_{C^1}$ is controlled by $\|\lambda_{\Phi}\|_{\infty}$ and hence controlled by $\|\lambda_{\Phi}'\|_{\infty}$. So we have
\begin{align*}
\|\nabla\Psi\|_\infty \leq C_0  \|\lambda_{\Phi}'\|_{\infty}.
\end{align*} 
Plugging this into \eqref{eq:dvPhimzeta} and taking the supremum over all test functions $\zeta\in C^{\infty}_c(V)$ with $\|\zeta\|_{\infty}\leq 1$ we deduce
\begin{align*}
|\nabla\cdot\Phi(m)|(V)\leq C_0 \,\liminf_{n\to\infty} I_{\e_n}(m_n;V)\quad\forall \Phi\in\mathrm{ENT}\text{ with }\|\lambda'_\Phi\|_{\infty}\leq 1.
\end{align*}
Hence for any open subset $U\subset\Omega$ and any disjoint open subsets $V_1,\ldots, V_k\subset U$ and entropies $\Phi_1,\ldots,\Phi_k\in \mathrm{ENT}$ with $\|\lambda'_{\Phi_j}\|_{\infty}\leq 1$ we have
\begin{align*}
\sum_j |\nabla\cdot\Phi_j(m)|(V_j) &\leq C_0 \sum_j \liminf_{n\to\infty} I_{\e_n}(m_n;V_j)\\
&\leq C_0  \liminf_{n\to\infty} I_{\e_n}(m_n;U).
\end{align*}
Given any disjoint compact sets $A_1,\ldots,A_k\subset U$, we can find disjoints open sets containing them, and so for entropies $\Phi_1,\ldots,\Phi_k\in \mathrm{ENT}$ with $\|\lambda'_{\Phi_j}\|_{\infty}\leq 1$ we have
\begin{align*}
\sum_j |\nabla\cdot\Phi_j(m)|(A_j) 
&\leq C_0  \liminf_{n\to\infty} I_{\e_n}(m_n;U).
\end{align*}
By inner regularity of the Radon measures $\nabla\cdot\Phi(m)$ this is in fact valid for any disjoint measurable sets $A_1,\ldots,A_k$, and then for any countable disjoint family of measurable sets $\lbrace A_j\rbrace$.
Recalling the definition
\cite[Definition~1.68]{ambrosio} of the lowest upper bound measure, this implies the lower bound \eqref{eq:lowerbound}.
\end{proof}

\section{Compactness}\label{s:compactness}

In this section we prove Theorem~\ref{t:comp}: let $\{m_n\}\subset W^{1,2}(\Omega;\R^2)$ satisfy $\nabla\cdot m_n=0$ and $\sup_n I_{\e_n}(m_n)<\infty$, then $\{m_n\}$ is precompact in $L^2(\Omega)$. The proof follows very closely the arguments in \cite[Proposition~1.2]{mul2}. We only briefly sketch the main ideas and highlight the steps that require adaptation. We refer to \cite{mul2} for the details that stay unchanged. 

The proof consists in showing that any Young measure $\lbrace\mu_x\rbrace_{x\in\Omega}$ generated by a subsequence of $\{m_n\}$ must be a family of Dirac measures. As in \cite[(3.17)]{mul2}, the energy bound $\sup_n I_{\e_n}(m_n)<\infty$ implies that $\mu_x$ is concentrated on $\partial\Bs$ for a.e. $x\in\Omega$.  The first main step is to prove that, for any entropy $\Phi\in\mathrm{ENT}$   the sequence $\nabla\cdot\widehat\Phi(m_n)$ is precompact in $H^{-1}(\Omega)$. This follows from the identity obtained in Lemma~\ref{l:extPhi}, exactly as in \cite[(3.1)]{mul2}.
The div-curl lemma therefore implies that for any entropies $\Phi_1,\Phi_2\in\mathrm{ENT}$, the weak* limit of the product $\widehat\Phi_1(m_n)\cdot i\widehat\Phi_2(m_n)$ in measures is the product of the weak limits of $\widehat\Phi_1(m_n)$ and $i\widehat\Phi_2(m_n)$ in $L^2(\Omega)$.  Hence, for a.e. $x\in\Omega$,   $\mu=\mu_x$ is a probability measure concentrated on $\partial\Bs$ and satisfying
\begin{align}\label{eq:Phi1Phi2mu}
\int \Phi_1\cdot i \Phi_2\, d\mu = \int \Phi_1\,  d\mu \cdot i  \int \Phi_2\,  d\mu.
\end{align}
The conclusion of Theorem~\ref{t:comp} then follows from the next Lemma, which is the counterpart of \cite[Lemma~2.6]{mul2}.
\begin{lem}
	\label{l:muDirac}
	Let $\mu$ be a probability measure on $\R^2$ that is supported on $\partial\Bs$ and satisfies \eqref{eq:Phi1Phi2mu}
	for all $\Phi_1, \Phi_2\in\mathrm{ENT}$. Then $\mu$ is a Dirac measure.
\end{lem}

\begin{rem} 
In the euclidean case $\|\cdot\|=|\cdot|$, building on earlier work by Aviles and Giga \cite{avilesgig},  Jin and Kohn \cite{kohn} introduced two fundamental entropies $\Sigma_1, \Sigma_2:\mathbb{S}^1\to\R^2$ given by 
\begin{align*}
\Sigma_1\lt(e^{i\theta}\rt)
= \frac i 2\lt(\frac{e^{i 3\theta}}{3}+e^{-i\theta}\rt),
\quad \Sigma_2\lt(e^{i\theta}\rt)=\frac 12\lt(\frac{e^{i3\theta}}{3}-e^{-i\theta}\rt).
\end{align*}
The proof of compactness given by \cite{ADM} uses only $\Sigma_1,\Sigma_2$ but is somewhat intricate, and the one in \cite{mul2} uses an infinite family of entropies.
We indicate here a somewhat shorter proof using only $\Sigma_1,\Sigma_2$  and  \v{S}ver\'ak's theorem \cite{sverak}.
To see this, rewrite \eqref{eq:Phi1Phi2mu} applied to $\Sigma_1,\Sigma_2$ as
\begin{align*}
\det\left(\int_{\R^{2\times 2}} X \,d\nu(X)\right) = \int_{\R^{2\times 2}} \det (X)\, d\nu(X),
\end{align*} 
where $\nu=\PI_{\sharp}\mu$ is the pushforward of $\mu$ by the matrix-valued map $\PI\colon\mathbb{S}^1\rightarrow \R^{2\times 2}$ whose rows are $\Sigma_1,\Sigma_2$.
Hence $\nu$ is a Null Lagrangian measure  (in the sense of 
\cite{LPnull}) supported on $K=\PI(\mathbb S^1)$ (which is the same set as in \cite[(42)]{LP}).
 By \cite[Lemma~7]{LP}, the set $K$ has no Rank-$1$ connections, so \cite[Lemma~3]{sverak} ensures that $\nu$ is a Dirac measure.
\end{rem}

The proof of Lemma~\ref{l:muDirac} follows  closely \cite[Lemma~2.6]{mul2}. Nevertheless we give some details, because this is where the crucial assumption that $\Bs$ is strictly convex is used. In the proof we will need the following construction.
\begin{lem}
	\label{l:genent}
	Given $\xi=\gamma(\theta_0)\in \partial \Bs$, define $\Phi^{\xi}:\partial\Bs\to\R^2$ by
\begin{align*}
\Phi^{\xi}(z)=\mathbf 1_{z\cdot i\gamma(\theta_0)>0}\,\gamma'(\theta_0)=\mathbf 1_{z\cdot i\xi>0}\,in_\Bs(\xi),
\end{align*}
where $n_{ \Bs}(\xi)$ denotes the outer unit normal to $\partial \Bs$ at $\xi$. 
	Then $\Phi^{\xi}$ is a generalized entropy for the equation \eqref{eq:geneikon} in the sense that there exists a sequence $\lt\{\Phi^{\xi}_{\delta}\rt\}_{\delta>0}\subset\mathrm{ENT}$ that is uniformly bounded and satisfies
	\begin{equation}
		\label{eqb2}
		\Phi^{\xi}_{\delta}(z)\rightarrow\Phi^{\xi}(z)\text{ for all }z\in \partial\Bs.
	\end{equation}
\end{lem}

\begin{proof}
For any $\lambda\in C^1(\R/2\pi\Z)$ such that $\int_{\R/2\pi\Z}\lambda(\theta)\gamma'(\theta)\, d\theta =0$, the map $\Phi\colon\partial \Bs\to\R^2$ given by
\begin{align*}
\Phi(\gamma(\theta))=\int_{\theta_0}^\theta \lambda(t)\gamma'(t)\, dt,
\end{align*}
is well defined and belongs to $\mathrm{ENT}$. We define a sequence of functions $\lambda_\delta$ such that the corresponding $\Phi=\Phi^\xi_\delta$ has the desired properties. We fix a smooth nonnegative kernel $\rho\in C_c^\infty(\R)$ with support $\supp\rho\subset (0,1)$ and unit integral $\int\rho=1$, denote $\rho_\delta(t)=\delta^{-1}\rho(t/\delta)$, and define a function $\hat\lambda_\delta\in C^\infty(\R/2\pi\Z)$ by setting
\begin{align*}
\hat\lambda_\delta(\theta)=  \rho_\delta(\theta-\theta_0) 
+ \rho_\delta(\pi+\theta_0-\theta) 
\qquad\text{for }\theta_0 <  \theta \leq \theta_0 +2\pi,
\end{align*}
and $\hat\lambda_\delta$ extended as a $2\pi$-periodic function. Note that $\hat\lambda_\delta$ is supported in 
\begin{align*}
(\theta_0,\theta_0+\delta) \cup (\theta_0+\pi-\delta,\theta_0 +\pi) +2\pi\Z.
\end{align*}
Moreover, the map $\Psi_\delta\colon (\theta_0,\theta_0+2\pi]\to\R^2$ defined by
\begin{align*}
\Psi_\delta(\theta)=\int_{\theta_0}^\theta \hat\lambda_\delta(t)\gamma'(t)\, dt\qquad\text{for }\theta_0 <  \theta \leq \theta_0 +2\pi,
\end{align*}
satisfies $|\Psi_\delta|\leq 2$ since $|\gamma'|=1$, and
\begin{align*}
\Psi_\delta(\theta)
& 
=
\begin{cases}
\int_\R \rho_\delta(t-\theta_0)\gamma'(t)\, dt &\text{ if }\theta_0+\delta <\theta <\theta_0+\pi-\delta,\\
\int_\R \rho_\delta(t-\theta_0)\gamma'(t)\, dt
+\int_\R \rho_\delta(\theta_0+\pi-t)\gamma'(t)\, dt
&\text{ if }\theta_0+\pi \leq \theta \leq  \theta_0+2\pi.\\
\end{cases}
\end{align*}
Using that $\gamma'$ is continuous and $\gamma'(\theta_0+\pi)=-\gamma'(\theta_0)$ we obtain, for any $\theta\in (\theta_0,\theta_0+2\pi]$, the limit
\begin{align*}
\lim_{\delta\to 0}\Psi_\delta(\theta)
&=
\begin{cases}
\gamma'(\theta_0)  &\quad\text{ if }\theta_0 <\theta <\theta_0+\pi ,\\
0
&\quad \text{ if }\theta_0+\pi \leq \theta \leq \theta_0+2\pi.\\
\end{cases}
\end{align*}
This corresponds exactly to $\Phi^\xi(\gamma(\theta))$. The only issue is that $\hat\lambda_\delta$ does not satisfy the constraint 
$\int_{\R/2\pi\Z}\hat\lambda_\delta(\theta)\gamma'(\theta)\, d\theta =0$, hence $\Psi_\delta$ cannot be extended to a $2\pi$-periodic function and does not define an entropy. So we need to modify $\hat\lambda_\delta$.\footnote{Note in order to get the pointwise convergence in \eqref{eqb2} for every $z\in \partial \Bs$ we can not define $\hat\lm_{\delta}$ via the standard symmetric (across zero) kernel centered on $\theta_0$, $\theta_0+\pi$. This is why we do the two step procedure of defining $\hat\lambda_\delta$ then modifying it.} We claim that there exists $\mu_\delta\in C^1(\R/2\pi\Z)$ such that
\begin{align}\label{eq:mudelta}
\begin{aligned}
&\int_{\R/2\pi\Z}\mu_\delta(\theta)\gamma'(\theta)\, d\theta
=
\int_{\R/2\pi\Z}\hat\lambda_\delta(\theta)\gamma'(\theta)\, d\theta,\\
\text{and }&\int_{\R/2\pi\Z}|\mu_\delta(\theta)|\,d\theta \longrightarrow 0\qquad\text{as }\delta\to 0. 
\end{aligned}
\end{align}
Granted \eqref{eq:mudelta}, we define $\lambda_\delta =\hat\lambda_\delta(\theta)-\mu_\delta$. This function does satisfy $\int_{\R/2\pi\Z}\lambda_\delta(\theta)\gamma'(\theta)\, d\theta =0$ thanks to \eqref{eq:mudelta}, so  the formula
\begin{align*}
\Phi_\delta^\xi(\gamma(\theta))=\int_{\theta_0}^\theta \lambda(t)\gamma'(t)\, dt,
\end{align*}
defines an entropy $\Phi_\delta^\xi\in\mathrm{ENT}$. Moreover for any $\theta\in (\theta_0,\theta_0+2\pi]$ we have
\begin{align*}
|\Psi_\delta(\theta)-\Phi_\delta^\xi(\gamma(\theta))|\leq 
\int_{\R/2\pi\Z}|\mu_\delta(\theta)|\,d\theta \longrightarrow 0\qquad\text{as }\delta\to 0. 
\end{align*}
Thanks to the convergence of $\Psi_\delta$ established above, this implies $\Phi_\delta^\xi(z)\to \Phi^\xi(z)$ for all $z\in\partial\Bs$, and uniform boundedness of $\Phi_\delta^\xi$ follows from the uniform boundedness of $\Psi_\delta$. Therefore the proof of Lemma~\ref{l:genent} will be complete once we prove the existence of  $\mu_\delta\in C^1(\R/2\pi\Z)$  satisfying \eqref{eq:mudelta}. 

This construction is possible thanks to the fact that
\begin{align*}
v^\delta:=\int_{\R/2\pi\Z}\hat\lambda_\delta(\theta)\gamma'(\theta)\, d\theta\longrightarrow \gamma'(\theta_0)+\gamma'(\theta_0+\pi)=0,
\end{align*}
as $\delta\to 0$.
To explicitly construct $\mu_\delta$,  we introduce a (small) parameter $\eta>0$, to be fixed later, and two functions $f_1^\eta,f_2^\eta\in C^1(\R/2\pi\Z)$ such that
\begin{align*}
\|f_j^\eta-\gamma_j'\|_{L^2(\R/2\pi\Z)}\leq\eta\qquad\text{for }j=1,2,
\end{align*}
and look for $\mu_\delta$ in the form
\begin{align*}
\mu_\delta (\theta)=\alpha_\delta f_1^\eta(\theta) +\beta_\delta f_2^\eta(\theta),\qquad \alpha_\delta,\beta_\delta\in\R.
\end{align*}
With these notations, the constraint 
$\int \mu_\delta \gamma' 
=
\int \hat\lambda_\delta \gamma' $ in \eqref{eq:mudelta} turns into
\begin{align*}
v^\delta=A_\eta \left(\begin{array}{c} \alpha_\delta\\ \beta_\delta\end{array}\right),
\qquad
A_\eta=\left(\begin{array}{cc}
\int f_1^\eta\gamma'_1 & \int f_2^\eta\gamma'_1\\
\int f_1^\eta\gamma'_2 & \int f_2^\eta\gamma'_2
\end{array}
\right).
\end{align*}
Next we show that we may fix $\eta>0$ such that $A_\eta$ is invertible. As a consequence, defining $(\alpha_\delta,\beta_\delta)^{T}=A_\eta^{-1}v^\delta$ ensures that $\mu_\delta$ satisfies the constraint 
$\int \mu_\delta \gamma' 
=
\int \hat\lambda_\delta \gamma' $
 in \eqref{eq:mudelta}, and the convergence $\mu_\delta\to 0$ in $L^1$ follows from $v^\delta\to 0$. This concludes the proof of \eqref{eq:mudelta} and of Lemma~\ref{l:genent}.

It remains to prove that $A_\eta$ is invertible for small enough $\eta>0$. To that end we remark that
thanks to the convergence $f_j^\eta\to \gamma_j'$ in $L^2$, we have
\begin{align*}
A_\eta \longrightarrow A_0 =
\left(\begin{array}{cc}
\int (\gamma'_1)^2 & \int \gamma'_2\gamma'_1\\
\int \gamma'_1\gamma'_2 & \int (\gamma'_2)^2
\end{array}
\right)\qquad\text{ as }\eta\to 0.
\end{align*}
The Cauchy-Schwarz inequality $(\int \gamma'_1\gamma'_2)^2\leq\int (\gamma'_1)^2\int (\gamma'_2)^2$ ensures that $\det(A_0)\geq 0$, and in fact $\det(A_0)>0$ because equality cannot occur in the Cauchy-Schwarz inequality: otherwise $\gamma'_1,\gamma'_2$ would be colinear in $L^2$, implying that $\gamma'$ takes values in a fixed line, which is incompatible with it being the unit tangent of $\partial\Bs$. So the matrix $A_0$ is invertible, and we may fix $\eta>0$ such that $A_\eta$ is invertible. 
\end{proof}

With the construction of Lemma~\ref{l:genent} at hand, we turn to the proof of Lemma~\ref{l:muDirac}.

\begin{proof}[Proof of Lemma~\ref{l:muDirac}]
	Let $\xi_1,\xi_2\in \partial \Bs$.
	By Lemma \ref{l:genent}, for $j=1,2$,  we can find $\lt\{\Phi_{\delta}^{\xi_j}\rt\}_{\delta}\subset \mathrm{ENT}$ such that 
	\begin{equation*}
		\Phi^{\xi_j}_{\delta}(z)\overset{\delta\rightarrow 0}{\rightarrow} \Phi^{\xi_j}(z)\qd\text{ for all }z\in \partial \Bs.
	\end{equation*}
 Also recall that $\lt\{\Phi^{\xi_j}_{\delta}\rt\}$ is uniformly bounded. Hence, applying \eqref{eq:Phi1Phi2mu} to $\Phi_j=\Phi^{\xi_j}_\delta$ and passing to the limit $\delta\to 0$ we obtain, by dominated convergence,
	\begin{equation*}
		\int_{\partial \Bs}  \Phi^{\xi_1}(z)\cdot i  \Phi^{\xi_2}(z)  d\mu(z) = \lt(\int_{\partial \Bs}  \Phi^{\xi_1}(z) d\mu(z)\rt)\cdot i\lt(\int_{\partial \Bs}  \Phi^{\xi_2}(z)  d\mu(z)\rt).
	\end{equation*}
In other words, recalling the definition of $\Phi^\xi$ in Lemma~\ref{l:genent},  for any $\xi_1,\xi_2\in \partial \Bs$ we have
	\begin{align*}
& in_\Bs(\xi_1)\cdot n_\Bs(\xi_2) \,\mu( \lbrace z\cdot i\xi_1>0\rbrace\cap \lbrace z\cdot i\xi_2>0\rbrace)
\nonumber\\
&
= in_\Bs(\xi_1)\cdot n_\Bs(\xi_2) \,\mu( \lbrace z\cdot i\xi_1>0\rbrace)\,\mu( \lbrace z\cdot i\xi_2>0\rbrace).
	\end{align*}
	(Here and in the rest of the proof, we use the shortened notation $\lbrace z\cdot i\xi>0\rbrace$ to denote the subset of points $z\in\partial \Bs$ satisfying this inequality.)
	By strict convexity of $\partial\Bs$, for  $\xi_1\ne \pm \xi_2$ we have $in_{\Bs}(\xi_1)\cdot n_{\Bs}(\xi_2)\ne 0$, and the last equation becomes
	\begin{align*}
\mu( \lbrace z\cdot i\xi_1>0\rbrace\cap \lbrace z\cdot i\xi_2>0\rbrace)
= \mu( \lbrace z\cdot i\xi_1>0\rbrace)\,\mu( \lbrace z\cdot i\xi_2>0\rbrace).
	\end{align*}
From that point on the proof follows exactly \cite[Lemma~2.6]{mul2}, for the reader's convenience we recall here the short argument. Letting $\xi_2\to\xi_1$ with $\xi_2\neq \pm\xi_1$   
we obtain
\begin{equation*}
		\mu\lt( \lbrace z\cdot i\xi_1>0\rbrace\rt)\leq\mu\lt(   \lbrace z\cdot i\xi_1>0\rbrace\rt)\mu\lt(   \lbrace z\cdot i\xi_1\geq0\rbrace\rt),
\end{equation*}
which implies
\begin{equation*}
		\mu\lt( \lbrace z\cdot i\xi>0\rbrace\rt)=0\text{ or }\mu\lt( \lbrace z\cdot i\xi\geq0\rbrace\rt)=1\text{ for all }\xi\in \partial \Bs.
\end{equation*}
This is equivalent to
\begin{equation}
	\label{eqb26}
	\mu\lt( \lbrace z\cdot i\xi>0\rbrace\rt)=0\text{ or }\mu\lt( \lbrace z\cdot i\xi<0\rbrace\rt)=0\text{ for all }\xi\in \partial \Bs.
\end{equation}
This implies that $\mu$ is a Dirac measure: otherwise we can find $\xi\in \partial \Bs$ such that 
	\begin{equation*}
		\mu\lt(\lbrace z\cdot i\xi>0\rbrace\rt)>0 \text{ and }\mu\lt(\lbrace z\cdot i\xi<0\rbrace\rt)>0,
	\end{equation*}
which contradicts \eqref{eqb26}.
\end{proof}

%
%
%
%

%

\section{Zero-energy states}\label{s:zeroenergy}

In this section we prove Theorem~\ref{T:zeroenergy} on zero-energy states: 
let $m\colon\Omega\to \R^2$ satisfy $m=\lim_{n\to\infty} m_n$ in $L^2(\Omega)$ for some $\{m_n\}\subset W^{1,2}(\Omega;\R^2)$ 
with $\na\cdot m_n=0$ and $\lim_{n\to \infty} I_{\e_n}(m_n)=0$ for some $\e_n\to 0$, 
then $m$ must be continuous outside a locally finite set of vortices associated to the norm $\|\cdot\|$.

The first step, similar to \cite[Proposition~1.1]{otto} is to obtain the kinetic formulation
\begin{align}\label{eq:kinzero}
\gamma'(t)\cdot \na_x\mathbf 1_{m(x)\cdot i\gamma(t)>0} = 0\quad\text{ in }\mathcal{D}'(\Omega)\text{ for all }t\in\R.
\end{align}
This follows from the fact that entropy productions vanish: thanks to Proposition~\ref{p:lowerbound}, $\dv\Phi(m)=0$ in $\mathcal{D}'(\Omega)$ for all $\Phi\in\mathrm{ENT}$. For any $t\in\R$, we may apply this to the entropies $\Phi_\delta^{\gamma(t)}$ provided by Lemma~\ref{l:genent}, hence
\begin{align*}
\int_\Omega \Phi_\delta^{\gamma(t)}(m(x))\cdot \nabla\zeta(x)\, dx=0\qquad\text{for all }\zeta\in C_c^\infty(\Omega).
\end{align*}
Thanks to the pointwise convergence  $\Phi_\delta^{\gamma(t)}(z)\to \Phi^{\gamma(t)}(z)=
\mathbf 1_{z\cdot i\gamma(t)>0}\gamma'(t)$ and the uniform boundedness of $\Phi_\delta^{\gamma(t)}$, we can pass to the limit $\delta\to 0$ by dominated convergence, and deduce
\begin{align*}
\int_\Omega \mathbf 1_{m(x)\cdot i\gamma(t)>0}\,\gamma'(t)\cdot \nabla\zeta(x)\, dx=0\qquad\text{for all }\zeta\in C_c^\infty(\Omega),
\end{align*}
which is exactly \eqref{eq:kinzero}.

Next we define $\tilde m\colon\Omega\to\mathbb S^1$ by setting
\begin{align*}
\tilde m =n_\Bs(m),\quad\text{where } n_\Bs\colon\partial\Bs\to\mathbb S^1\text{ is the outer unit normal to }\partial\Bs.
\end{align*}
The symmetries of $\Bs$ ensure that, for 
 any $z\in\partial \Bs$,
\begin{align*}
z\cdot i\gamma(t) >0 \qquad\Longleftrightarrow\qquad n_\Bs(z)\cdot \gamma'(t)>0.
\end{align*}
Therefore, for a fixed $t\in \R$, and $\theta_t=\alpha(t)$ where $\alpha\colon\R\to\R$ is the (unique up to an additive constant) continuous function such that $\gamma'(t)=e^{i\alpha(t)}$, we have
\begin{align*}
\mathbf 1_{m(x)\cdot i\gamma(t)>0}\,\gamma'(t) 
= \mathbf 1_{\tilde m(x)\cdot e^{i\theta_t}>0}\, e^{i\theta_t}.
\end{align*}
As $t\mapsto\theta_t$ is a bijection from $\R$ into itself we deduce from \eqref{eq:kinzero} that $\tilde m$ solves the kinetic equation
\begin{align*}
e^{i\theta}\cdot\nabla_x\mathbf 1_{\tilde m(x)\cdot e^{i\theta}>0} =0
\quad\text{ in }\mathcal{D}'(\Omega)\text{ for all } \theta\in\R.
\end{align*}
This is the kinetic formulation that characterizes zero-energy states for the classical Aviles-Giga functional: it follows from \cite[Theorem~1.3]{otto} that $\tilde m$ is  locally Lipschitz outside a locally finite set. Moreover, in any convex neighborhood of a singularity $x_0$ we have $\tilde m(x)=\beta i(x-x_0)/|x-x_0|$ for some $\beta\in\lbrace \pm 1\rbrace$. 

Note that, since $\partial \Bs$ is $C^1$ and strictly convex, the map
 $n_\Bs\colon\partial\Bs\to\mathbb S^1$ is a homeomorphism. We deduce that $m=n_\Bs^{-1}(\tilde m)$  is continuous outside a locally finite set. Moreover, from \cite[Proposition~2.2]{boch} we know that $n_\Bs^{-1}(x/|x|)=V_\Bs(x)$ for any $x\in\R^2$, where $V_\Bs=\nabla\|\cdot\|_*$ is the vortex associated to $\|\cdot\|$ and $\|\cdot \|_{*}$ is the dual norm of $\|\cdot \|$, and $V_\Bs(-x)=-V_\Bs(x)$, so we deduce $m(x)=\beta V_\Bs(i(x-x_0))$ in any convex neighborhood of a singularity $x_0$. This concludes the proof of Theorem~\ref{T:zeroenergy}. 
 
To prove the assertion of Remark~\ref{r:zeroenergytypep}, simply note that $n_\Bs^{-1}$ is $1/(p-1)$-H\"older whenever $\|\cdot\|$ is of power type $p$  for $p\geq 2$ \cite[Theorem~2.6]{boch}.


\section{Regularity estimates}\label{s:reg}

In this section we give the proofs of Theorems \ref{thmbes} and \ref{thmbesrev} in Subsections \ref{ss:reg} and \ref{ss:revreg}, respectively. The proof of Theorem \ref{t:regconvBV} relies on explicit calculations in the $BV$ setting and is postponed to Subsection \ref{besbv}.

\subsection{Finite entropy production implies regularity estimates}
\label{ss:reg}

In this subsection we find it more convenient to  work with maps $m\colon\Omega\to\mathbb{R}^2$ solving
\begin{align}\label{eq:mX}
	|m|=1\text{ a.e.},\quad \nabla\cdot X(m)=0\text{ in }\mathcal D'(\Omega).
\end{align}
Solutions of \eqref{eq:mX} and of the generalized Eikonal equation \eqref{eq:geneikon} are in correspondence via the Lipschitz homeomorphism $X$ defined in \eqref{eq:X}. Specifically, a map $m$ solves \eqref{eq:mX} if and only if $\overline m = X(m)$ solves \eqref{eq:geneikon}.

This transformation also induces  a correspondence between entropies. 
	A $C^1$ map $\Phi\colon \mathbb S^1\to\R^2$ is an entropy for equation \eqref{eq:mX} if and only if
	$\nabla\cdot\Phi(m)=0$ for any smooth solution of \eqref{eq:mX}. 
	It is an exercise to see that this is equivalent to $(d/d\theta)\Phi(e^{i\theta})$ being colinear to $\gamma'(\theta)$. 
	As in the definition of $\mathrm{ENT}$ in \eqref{eq:ENT}, we consider the subclass of entropies where we require a $C^1$ colinearity coefficient:
	\begin{align*}
		\frac{d}{d\theta}\Phi(e^{i\theta})=\lambda(\theta)\gamma'(\theta)\quad\text{for some }\lambda\in C^1(\R/2\pi\Z).
	\end{align*}
	Therefore, $\Phi$ is an entropy for equation \eqref{eq:mX}  in this subclass
	if and only if $\overline\Phi=\Phi\circ X^{-1}\colon\partial\Bs\to\R^2$  belongs to $\mathrm{ENT}$, as  follows directly from the definition \eqref{eq:ENT} of the class $\mathrm{ENT}$. Moreover  we have $\overline\Phi(\overline m)=\Phi(m)$. In particular, the entropy productions $\nabla\cdot\Phi(m)$ of a solution $m$ of \eqref{eq:mX} are measures if and only if the entropy productions $\nabla\cdot\overline\Phi(\overline m)$ of the solution $\overline m=X(m)$ of \eqref{eq:geneikon} are measures. Thanks to the above discussion, all results we prove in this section about solutions of \eqref{eq:mX} directly translate into corresponding results about solutions of \eqref{eq:geneikon}.

We use the family of entropies $\Phi_\psi$ for equation \eqref{eq:mX} given by
\begin{align}\label{eq:Phipsi}
	\Phi_\psi(e^{i\theta})=\int_{\mathbb{R}/ 2\pi \mathbb{Z}} \mathbf 1_{e^{i\theta}\cdot e^{is}>0}\,\psi(s)\gamma'(s-\pi/2)\, ds,\qquad\psi\in C^1(\R/2\pi\mathbb Z).
\end{align}
Note that since $\gamma'(t+\pi)=-\gamma'(t)$ we have that $(d/d\theta)\Phi_{\psi}\lt(e^{i\theta}\rt)=\lambda(\theta)\gamma'(\theta)$ with
$\lambda(\theta)=\psi(\theta+\pi/2)+\psi(\theta-\pi/2)$. 
Therefore \eqref{eq:Phipsi}  does define an entropy for equation \eqref{eq:mX}.

Recall that $\alpha\in C^0(\R)$ is  such that
\begin{align}\label{eq:alpha}
	e^{i\alpha(\theta)}=\gamma'(\theta)\qquad  \forall \theta\in\R.
\end{align}
The continuous function $\alpha$ is uniquely determined up to a constant, is  strictly increasing, and satisfies $\alpha(\theta+\pi)=\alpha(\theta)+\pi$ for all $\theta\in\R$.

\begin{prop}\label{p:besov}
	If $m$ satisfies \eqref{eq:mX} and
	\begin{align*}
		\dv\Phi_\psi(m)\in\mathcal M_{loc}(\Omega)\qquad\forall\psi\in C^1(\mathbb R/2\pi\mathbb Z),
	\end{align*}
	then we have
	\begin{align}
	\label{eqbes22}
		\sup_{0<|h|<\dist(\Omega',\partial\Omega)} \frac{1}{|h|}\int_{\Omega'} \Lambda(m(x),m(x+h))\, dx   
		<\infty\qquad\forall\Omega'\subset\subset\Omega,
	\end{align}
	where $\Lambda\colon\mathbb S^1\times \mathbb S^1\to [0,+\infty)$ is given by 
	\begin{align}\label{eq:Lambda1}
		\Lambda(e^{i\theta_1},e^{i\theta_2}) =\int_{\theta_1}^{\theta_2} \int_{\theta_1}^{\theta_2} |\alpha(t)-\alpha(s)|\, dt ds \quad\text{ for }|\theta_1-\theta_2|\leq \pi.
	\end{align}
	Moreover, for $m_1, m_2\in\mathbb{S}^1$, we have 
	\begin{align}
	\label{eqbes21}
		\Lambda(m_1,m_2)\gtrsim \delta^2\omega^{-1}(\delta/2),\qquad\delta=|m_1-m_2|,
	\end{align}
	where $\omega(\delta)=\sup \lbrace |\alpha^{-1}(t)-\alpha^{-1}(s)|\colon |t-s|<\delta \rbrace$ is the minimal modulus of continuity of $\alpha^{-1}$.
\end{prop}
Theorem \ref{thmbes} is a direct consequence of Proposition \ref{p:besov} and the correspondence between the generalized Eikonal equation \eqref{eq:geneikon} and equation \eqref{eq:mX}. In particular, the regularity estimate \eqref{eqbes2} of Theorem \ref{thmbes} is equivalent to \eqref{eqbes22} by noting that the function $\Pi$ defined in \eqref{eq:Lambda} satisfies $\Pi(\gamma(\theta_1),\gamma(\theta_2)) = \Lambda(e^{i\theta_1},e^{i\theta_2})$.
\begin{rem}
\label{remppower}
	If $\|\cdot\|$ is  in addition of power type $p$ for some $p\geq 2$ then $\alpha^{-1}$ is $1/(p-1)$-H\"older \cite[Theorem 2.6]{boch}, so $\omega^{-1}(\delta)\gtrsim \delta^{p-1}$ and therefore \eqref{eqbes22} and \eqref{eqbes21} imply
	\begin{align*}
		\sup_{0<|h|<\dist(\Omega',\partial\Omega)}\frac{1}{|h|}\int_{\Omega'} |m(x+h)-m(x)|^{p+1}\, dx  <\infty\qquad\forall\Omega'\subset\subset\Omega,
	\end{align*}
	that is, $m$ has the local Besov regularity $B^{\frac{1}{p+1}}_{p+1,\infty}$.
\end{rem}

\begin{rem}
\label{remlam}
	The function $\alpha$ is increasing and therefore $D\alpha$ is a positive measure. Hence for $0\leq \theta_2-\theta_1\leq\pi$, 
one can rewrite  the quantity $\Lambda$ defined in \eqref{eq:Lambda1} as
	\begin{align*}
		\Lambda(e^{i\theta_1},e^{i\theta_2})&= \iiint_{[\theta_1,\theta_2]^3}\left(\mathbf 1_{s<\tau<t} +\mathbf 1_{s>\tau>t}\right)\, dtds\, D\alpha(d\tau)
		\nonumber\\
		&=2\int_{\theta_1}^{\theta_2} (\tau-\theta_1)(\theta_2-\tau)\, D\alpha(d\tau)
		\nonumber\\
		&\geq \frac {2}{9} (\theta_2-\theta_1)^2 D\alpha([\theta_1+(\theta_2-\theta_1)/3,\theta_1+2(\theta_2-\theta_1)/3]).
	\end{align*}
(We used the identity $\iint_{[\theta_1,\theta_2]^2}\left(\mathbf 1_{s<\tau<t} +\mathbf 1_{s>\tau>t}\right)\, dtds=2(\tau-\theta_1)(\theta_2-\tau)$ to obtain the second equality.)
\end{rem}

\begin{proof}[Proof of Proposition~\ref{p:besov}]
	The proof is a direct combination of Lemmas~\ref{l:kin} to \ref{l:lowerLambda} below.
\end{proof}

\begin{lem}\label{l:kin}
	If $m$ satisfies \eqref{eq:mX} and
	\begin{align*}
		\dv\Phi_\psi(m)\in\mathcal M(\Omega)\qquad\forall\psi\in C^1\lt(\mathbb{R}/2\pi\mathbb{Z}\rt),
	\end{align*}
	then $m$ satisfies the kinetic equation
	\begin{align}\label{eq:kin}
		\gamma'(s-\pi/2)\cdot \nabla_x \mathbf 1_{m(x)\cdot e^{is}>0} =\partial_s\sigma(s,x)\qquad \text{in }\mathcal{D}'(\mathbb{R}/2\pi\mathbb{Z}\times\Omega),
	\end{align}
	for some $\sigma\in\mathcal M(\mathbb{R}/2\pi\mathbb{Z}\times\Omega)$.
\end{lem}
\begin{proof}
For any fixed $\zeta\in C^{\infty}_c(\Omega)$, 
the operator $T_{\zeta}:\psi \mapsto \la \na \cdot \Phi_{\psi}(m), \zeta\ra $ is a linear operator from $C^1(\R/2\pi\mathbb{Z})$ to $\R$.  
Further the estimate $\lt|\la \na \cdot \Phi_{\psi}(m), \zeta\ra \rt|\lesssim \|\psi\|_{C^1\lt( \mathbb{R}/2\pi\mathbb{Z}\rt)} \|\na \zeta\|_{L^\infty(\Omega)}$ 
implies that $T_{\zeta}$ is a bounded linear operator. 
Thus the same Banach-Steinhaus argument as in \cite[Lemma~3.4]{GL} provides the bound
	\begin{align*}
		|\langle \dv\Phi_\psi(m),\zeta\rangle | \lesssim \|\psi\|_{C^1(\mathbb{R}/2\pi\mathbb{Z})}\|\zeta\|_{L^\infty(\Omega)},
	\end{align*}
	for all $\zeta\in C^0_c(\Omega)$ and $\psi\in C^1(\mathbb{R}/2\pi\mathbb{Z})$. Moreover, when $\psi$ is a constant $\psi\equiv c$, we have $\dv\Phi_\psi(m)=2c\dv X(m)=0$, so in the above we can consider the quotient space $C^1(\mathbb{R}/2\pi\mathbb{Z})/\R\approx C^0(\R/2\pi\Z)$. Explicitly, for any $f\in C^0(\mathbb{R}/2\pi\mathbb{Z})$ consider the function $\psi[f]\in C^1(\mathbb{R}/2\pi\mathbb{Z})$ given by
	\begin{align*}
		\psi[f](t)=\int_0^t\left(f-\Xint{-}_{\mathbb{R}/2\pi\mathbb{Z}}f\right),
	\end{align*}
	then we have
	\begin{align*}
		\langle \dv\Phi_{\psi[f]}(m),\zeta\rangle \lesssim \|f\|_{L^\infty(\mathbb{R}/2\pi\mathbb{Z})}\|\zeta\|_{L^\infty(\Omega)},
	\end{align*}
	for all $\zeta\in C^0_c(\Omega)$ and $f\in C^0(\mathbb{R}/2\pi\mathbb{Z})$. As a consequence (see \cite[Appendix~B]{LPfacto} for a detailed proof) there exists a measure $\sigma\in \mathcal M(\mathbb{R}/2\pi\mathbb{Z}\times\Omega)$ such that, for all $\zeta\in C^{\infty}_c(\Omega)$ and $f\in C^0(\R/2\pi\mathbb{Z})$ with $\Xint{-}_{\mathbb{R}/2\pi\mathbb{Z}}f=0$, we have
	\begin{align*}
		\langle \dv\Phi_{\psi[f]}(m),\zeta\rangle=-\langle \sigma,f\otimes\zeta\rangle =\langle\partial_s\sigma,\psi[f]\otimes \zeta\rangle.
	\end{align*}
	From the the definition of $\Phi_\psi$ \eqref{eq:Phipsi}  we see that this is equivalent to
	\begin{align}
		\label{eqb11}
		\langle \gamma'(s-\pi/2) \cdot\nabla_x \mathbf 1_{m(x)\cdot e^{is}>0} -\partial_s\sigma,\psi(s)\zeta(x)\rangle =0,
	\end{align}
	for $\psi=\psi[f]$, that is, for all $\psi\in C^1(\mathbb{R}/2\pi\mathbb{Z})$ such that $\psi(0)=0$ and $\zeta\in C^{\infty}_c(\Omega)$. For constant $\psi =c$, equation  \eqref{eqb11} amounts to   $2 c\,\dv X(m)=0$, so it is in fact valid for any $\psi\in C^1(\mathbb{R}/2\pi\mathbb{Z})$. This proves the kinetic equation \eqref{eq:kin}.
\end{proof}

\begin{rem}
Note that the kinetic equation \eqref{eq:kin} only uniquely determines $\partial_s\sigma$. We may choose the unique $\sigma$ satisfying in addition $\langle\sigma(s,x), \zeta(x)\rangle=0$ for all $\zeta\in C^0_c(\Omega)$.
\end{rem}

\begin{lem}
	If $m$ satisfies $|m|=1$ a.e. and the kinetic equation \eqref{eq:kin}, then
	for any  $\varphi\in BV(\R/2\pi\Z)$ which is odd, i.e. $\varphi(-\theta)=-\varphi(\theta)$, the quantity
	\begin{align}\label{eq:Deltavarphi}
		\Delta_\varphi(e^{i\theta_1},e^{i\theta_2})&=
		\iint_{\R/2\pi\Z\times\R/2\pi\Z}\varphi(t-s)\,\gamma'(s-\pi/2)\wedge\gamma'(t-\pi/2)
		\\
		&\qquad\qquad 
		\left(\mathbf 1_{e^{is}\cdot e^{i\theta_2} >0}-\mathbf 1_{e^{is}\cdot e^{i\theta_1} >0}\right)\left(\mathbf 1_{e^{it}\cdot e^{i\theta_2} >0}-\mathbf 1_{e^{it}\cdot e^{i\theta_1} >0}\right)\, dtds, \nonumber
	\end{align}
	satisfies
	\begin{align*}
		\frac{1}{|h|}\int_{\Omega'}\Delta_\varphi(m(x),m(x+h))\, dx \lesssim \frac{\|\varphi\|_{L^1(\R/2\pi\mathbb{Z})}}{\dist(\Omega',\partial\Omega)}+|D\varphi|(\R/2\pi\mathbb{Z})
|\sigma|(\R/2\pi\mathbb{Z} \times\Omega),
	\end{align*}
	for all $\Omega'\subset\subset \Omega$ and $h\in\R^2$ such that $|h| <\dist(\Omega',\partial\Omega)$.
\end{lem}
\begin{proof}
	This essentially follows \cite[Lemma~3.9]{GL} in a slightly modified setting; we provide some details here for the reader's convenience. We set
	\begin{align*}
		\chi(t,x)=\mathbf 1_{e^{it}\cdot m(x)>0},
	\end{align*}
	and for a small parameter $\e>0$ we consider regularized (with respect to $x$) maps
	\begin{align*}
		\chi_\e=\chi *_x \rho_\e,\quad \sigma_\e =\sigma *_{x} \rho_\e,
	\end{align*}
	where $\rho_\e$ is a regularizing kernel. We have the regularized kinetic equation
	\begin{align*}
		\gamma'(s-\pi/2)\cdot\nabla_{x}\chi_\e=\partial_s\sigma_\e.
	\end{align*}
	
	Let $\Omega'\subset\subset\Omega$ be fixed and $h=u e$ for some $e\in\mathbb{S}^1$ and $u\in\R$ such that $|u|=|h|<\dist(\Omega',\partial\Omega)$. Without loss of generality, assume $e=e_1$. We denote by $\chi^u(t, x)=\chi(t, x+ue_1)$ and $D^u\chi(t,x)=\chi^u(t,x)-\chi(t,x)$.
	 Define the
	quantity
	\begin{align*}
		\Delta^{\e}_\varphi(x,u)&=\iint_{\R/2\pi\Z\times\R/2\pi\Z}\varphi(t-s)\,\gamma'(s-\pi/2)\wedge\gamma'(t-\pi/2) \\
		&\hspace{12em}
		D^{u}\chi_\e(s,x) D^{u}\chi_\e(t,x)\, dtds,
	\end{align*}
	for $x\in\Omega$ and $|u|+\epsilon < \dist(x,\partial\Omega)$.
	Note that, as $\e\to 0$, we have the pointwise limit
	\begin{align*}
		\Delta^{\e}_\varphi(x,u)\longrightarrow \Delta_\varphi(m(x),m(x+ ue_1))\quad\text{ for a.e. } x\in\Omega.
	\end{align*}
	A long but direct calculation (detailed in \cite[Lemma~3.9]{GL} in the case $\gamma(t)=e^{it}$) provides, for any smooth odd $\varphi$, the identity
	\begin{align*}
		\frac{\partial}{\partial u}\Delta^{\e}_\varphi & = I^{\ep}
		 +\nabla \cdot A^{\e},\\
		I^{\ep}& = 2\iint \varphi(t-s)\gamma'_2(t-\pi/2)\left[
		\chi_\e^{u}(t,x)\partial_s \sigma_\e(s,x)
		-\chi_\e(t,x)\partial_s\sigma_\e^{u}(s,x)
		\right] \, dtds,\\
		A_1^\e&=2\iint \varphi(t-s)\gamma'_2(t-\pi/2)\gamma'_1(s-\pi/2)\chi_\e^{u}(t,x)D^{u}\chi_\e(s,x)\, dtds, \\
		A_2^\e&= 2\iint \varphi(t-s)\gamma'_2(t-\pi/2)\gamma'_2(s-\pi/2)\chi_\e(t,x)\chi^u_\e(s,x)\, dtds.
	\end{align*}
	Note that $|A^{\e}|\lesssim \|\varphi\|_{L^1(\R/2\pi\mathbb{Z})}$. Integrating with respect to $u$ and against a smooth cut-off function in $x$  we deduce
	\begin{align*}
		\frac{1}{|h|}\int_{\Omega'} \Delta^{\e}_\varphi(x,u)\, dx
		&
		\lesssim
		\frac{\|\varphi\|_{L^1(\R/2\pi\mathbb{Z})}}{\dist(\Omega',\partial\Omega)} +
		|\sigma|(\R/2\pi\mathbb{Z}\times\Omega)\int_{\R/2\pi\mathbb{Z}} |\varphi'(t)|\, dt   .
	\end{align*}
	Letting $\e\to 0$ we infer
	\begin{align*}
		\frac{1}{|h|}\int_{\Omega'}\Delta_\varphi(m(x),m(x+u e_1))\, dx \lesssim \frac{\|\varphi\|_{L^1(\R/2\pi\mathbb{Z})}}{\dist(\Omega',\partial\Omega)}+
		|D\varphi|(\R/2\pi\mathbb{Z})|\sigma|(\R/2\pi\mathbb{Z}\times\Omega),
	\end{align*}
	for all smooth odd $\varphi$, and by approximation for any odd $\varphi\in BV(\R/2\pi\mathbb{Z})$. 
\end{proof}

\begin{lem}\label{l:lowerDelta}
	There exists an odd function $\varphi\in BV(\R/2\pi\Z)$  such that the quantity $\Delta_\varphi$ defined in \eqref{eq:Deltavarphi} satisfies
	$\Delta_\varphi\gtrsim \Lambda$, where $\Lambda$ is defined in \eqref{eq:Lambda1}.
\end{lem}
\begin{proof}
	We define an odd function $\varphi\in BV(\R/2\pi\Z)$ by setting
	\begin{align*}
		\varphi(\theta)=\begin{cases}
			1 &\text{ for }0< \theta<\delta,\\
			-1 &\text{ for }-\delta < \theta < 0,\\
			0 &\text{ for } \delta < |\theta | <\pi,
		\end{cases}
	\end{align*}
	where $\delta\in (0,\pi/2)$ is a parameter to be chosen later. 
	Recalling the definitions of $\Delta_\varphi$ \eqref{eq:Deltavarphi} and $\alpha$ \eqref{eq:alpha}, we have, for $m_1, m_2\in\mathbb{S}^1$,
	\begin{align}\label{eq:LambdaXi}
		\Delta_\varphi(m_1,m_2)&=\iint_{\dist_{\mathbb{S}^1}(e^{is},e^{it})<\delta} |\sin(\alpha(t-\pi/2)-\alpha(s-\pi/2))|\, \Xi(t)\,\Xi(s)\, dtds,\\
		\Xi(t)&=\Xi(t,m_1,m_2)=\mathbf 1_{e^{it}\cdot m_2>0}-\mathbf 1_{e^{it}\cdot m_1>0}.\nonumber
	\end{align}
	Here and in what follows we let $\dist_{\mathbb{S}^1}$ denote the geodesic distance in $\mathbb S^1$. The function $\Xi(\cdot,m_1,m_2)$ is supported in two opposite arcs of length $\dist_{\mathbb{S}^1}(m_1,m_2)$: 
	\begin{align*}
		\Xi(t,m_1,m_2)&=\mathbf 1_{t\in A} -\mathbf 1_{t\in -A}\qquad \text{for a.e. }t\in \R/2\pi\mathbb{Z}, \\
		A = A(m_1,m_2)&=\left\lbrace t\in \R/2\pi\mathbb{Z}\colon e^{it}\cdot m_2>0\text{ and }e^{it}\cdot m_1 < 0\right\rbrace.
	\end{align*}  
	If $\dist_{\mathbb{S}^1}(m_1,m_2)\leq  \pi-\delta$, then the distance between these arcs is at least $\delta$, so we have
	\begin{align*}
		\Xi(t)\Xi(s)&=\mathbf 1_{s,t\in A} + \mathbf 1_{s,t\in -A}\qquad\text{ for }\dist_{\mathbb{S}^1}(e^{is},e^{it})<\delta,
	\end{align*}
	and therefore 
	\begin{align*}
		\Delta_\varphi(m_1,m_2)&=\iint_{A\times A} \mathbf 1_{\dist_{\mathbb{S}^1}(e^{is},e^{it})<\delta} |\sin(\alpha(t-\pi/2)-\alpha(s-\pi/2))|\, dtds\\
		& \quad + \iint_{-A\times -A} \mathbf 1_{\dist_{\mathbb{S}^1}(e^{is},e^{it})<\delta} 
		|\sin(\alpha(t-\pi/2)-\alpha(s-\pi/2))|\, dtds\\
		&=2 \iint_{A\times A} \mathbf 1_{\dist_{\mathbb{S}^1}(e^{is},e^{it})<\delta} 
		|\sin(\alpha(t-\pi/2)-\alpha(s-\pi/2))|\, dtds.
	\end{align*}
	The last equality follows from $\alpha(t+\pi)=\alpha(t)+\pi$ for all $t\in \R/2\pi\mathbb{Z}$. Because $\alpha$ is uniformly continuous we may choose $\delta_0>0$ small enough to ensure that $|\alpha(t-\pi/2)-\alpha(s-\pi/2)|\leq\pi/2$ for $\dist_{\mathbb{S}^1}(e^{it},e^{is})<\delta$ provided $\delta\leq \delta_0$, and then we have
	\begin{align*}
		\Delta_\varphi(m_1,m_2)&\geq \frac 4\pi  \iint_{ A\times A} \mathbf 1_{\dist_{\mathbb{S}^1}(e^{is},e^{it})<\delta} 
		|\alpha(t-\pi/2)-\alpha(s-\pi/2)|\, dtds.
	\end{align*}
	Letting
	$m_1=e^{i\theta_1}$,
	$m_2=e^{i\theta_2}$ with $|\theta_1-\theta_2|=\dist_{\mathbb{S}^1}(m_1,m_2)$, this turns into
	\begin{align}\label{eq:lowDelta1}
		\Delta_\varphi(m_1,m_2)&\geq \frac 4\pi  \int_{\theta_1+\frac\pi 2}^{\theta_2 +\frac\pi 2}
		\int_{\theta_1+\frac\pi 2}^{\theta_2 +\frac\pi 2}  \mathbf 1_{|t-s|<\delta} |\alpha(t-\pi/2)-\alpha(s-\pi/2)|\, dtds \nonumber\\
		&=  \frac 4\pi  \int_{\theta_1}^{\theta_2}
		\int_{\theta_1}^{\theta_2}  \mathbf 1_{|t-s|<\delta} |\alpha(t)-\alpha(s)|\, dtds.
	\end{align}
	Recalling the definition \eqref{eq:Lambda1} of $\Lambda(m_1,m_2)$, we deduce
	\begin{align*}
		\Delta_\varphi(m_1,m_2)\geq
		\frac 4\pi \Lambda(m_1,m_2) \qquad \text{ if  }\dist_{\mathbb{S}^1}(m_1,m_2)\leq \delta.
	\end{align*}
	For $\delta<\dist_{\mathbb{S}^1}(m_1,m_2)\leq\pi-\delta$, from \eqref{eq:lowDelta1}  we have
	\begin{align*}
		\Delta_\varphi(m_1,m_2)& \geq C_1(\delta):=\frac 4\pi \inf_{\theta\in\R}\int_\theta^{\theta+\delta}\int_{\theta}^{\theta+\delta}|\alpha(t)-\alpha(s)|\, dtds >0
	\end{align*}
	(where this infimum is indeed positive because $\alpha$ is strictly increasing thanks to the strict convexity of $\Bs$),
	so
	\begin{align*}
		\Delta_{\varphi}(m_1,m_2)&\geq \frac{C_1(\delta)}{\sup_{\mathbb S^1\times\mathbb S^1}\Lambda} \Lambda(m_1,m_2)\qquad\text{ if }\delta<\dist_{\mathbb{S}^1}(m_1,m_2)\leq\pi-\delta.
	\end{align*}
	
	Finally we turn to the case $\pi- \delta <  \dist_{\mathbb{S}^1}(m_1,m_2) \leq\pi$, where the product $\Xi(s)\Xi(t)$ can take negative values. We have
	\begin{align*}
		\Xi(t)\Xi(s)&=\mathbf 1_{s,t\in A} + \mathbf 1_{s,t\in -A} -\mathbf 1_{t\in A,s\in -A} -\mathbf 1_{t\in -A,s\in A}\\
		&\geq \mathbf 1_{s,t\in A} + \mathbf 1_{s,t\in -A} -\mathbf 1_{s,t\in \hat A_\delta}-\mathbf 1_{s,t\in -\hat A_\delta},
	\end{align*}
	where $\hat A_\delta$ is the arc of length $2\delta$ given by
	\begin{align*}
		\hat A_\delta=\left\lbrace t\in\R/2\pi\mathbb{Z}\colon \dist_{\mathbb{S}^1}(e^{it},ie^{i\theta_0})<
		 \delta\right\rbrace,\qquad\theta_0=\frac{\theta_1+\theta_2}{2}-\frac\pi 2.
	\end{align*}
	Plugging this into the expression of $\Delta_\varphi$ in \eqref{eq:LambdaXi} and using again $\alpha(\theta+\pi)=\alpha(\theta)+\pi$ and $|\alpha(t-\pi/2)-\alpha(s-\pi/2)|\leq\pi/2$ for $\dist_{\mathbb{S}^1}(e^{is},e^{it})<\delta$, we find
	\begin{align}\label{eq:lowDeltalarge1}
		\Delta_\varphi(m_1,m_2)&\geq \frac{4}{\pi}\iint_{A\times A} \mathbf 1_{\dist_{\mathbb{S}^1}(e^{is},e^{it})<\delta} |\alpha(t-\pi/2)-\alpha(s-\pi/2)|\, dtds
		\nonumber \\
		&\quad - 2\iint_{\hat A_\delta\times\hat A_\delta}|\alpha(t-\pi/2)-\alpha(s-\pi/2)|\, dtds
		\nonumber \\
		& \geq  \frac{4}{\pi}\int_{\theta_0+\frac\delta 2}^{\theta_0+\pi -\frac\delta 2}\int_{\theta_0+\frac\delta 2}^{\theta_0+\pi -\frac\delta 2} \mathbf 1_{|t-s|<\delta} |\alpha(t)-\alpha(s)|\, dtds \nonumber \\
		&\quad - 2\int_{\theta_0 - \delta}^{\theta_0 +\delta}
		\int_{\theta_0 -\delta}^{\theta_0 +\delta}
		|\alpha(t)-\alpha(s)|\, dtds. 
	\end{align}
	Since $\alpha$ is increasing, its derivative $D\alpha$ is a nonnegative Radon measure. We use this to calculate
	\begin{align*}
		& \int_{\theta_0+\frac\delta 2}^{\theta_0+\pi -\frac\delta 2}\int_{\theta_0+\frac\delta 2}^{\theta_0+\pi -\frac\delta 2} \mathbf 1_{|t-s|<\delta} |\alpha(t)-\alpha(s)|\, dtds\\
		& =\int_{\theta_0+\frac\delta 2}^{\theta_0+\pi-\frac\delta 2} W(\tau) \, D\alpha(d\tau),
		\\
		W(\tau)&=\iint \mathbf 1_{|t-s|<\delta} \mathbf 1_{\theta_0+\delta/2 < t< \tau}
		\mathbf 1_{\tau < s< \theta_0+\pi-\delta/2}\,dtds\\
		&\quad + \iint \mathbf 1_{|t-s|<\delta}\mathbf 1_{\theta_0+\delta/2 < s< \tau}
			\mathbf 1_{\tau < t< \theta_0+\pi-\delta/2}
		\, dtds.
	\end{align*}
	For $\tau\in [\theta_0 +3\delta/2,\theta_0+\pi-3\delta/2]$ the quantity $W(\tau)$ is the sum of the areas of two isosceles right-angled triangles  of height $\delta$, so $W(\tau)=\delta^2$; see Figure \ref{tri}.
	\begin{figure}[t]
	\centering
	\includegraphics[width=.55\linewidth]{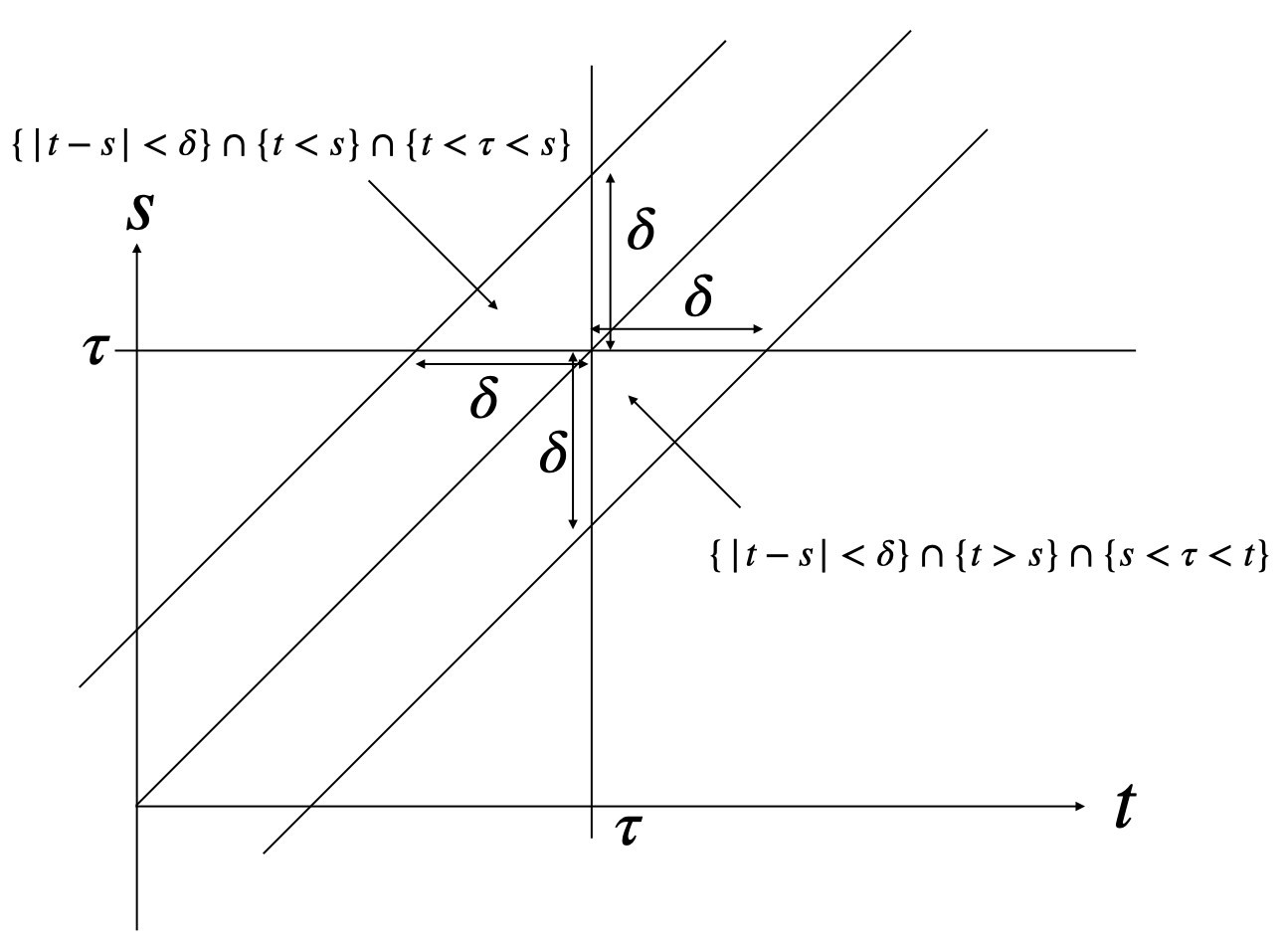}
	\caption{The two isosceles triangles in $W(\tau)$.}
	\label{tri}
	\end{figure}
	Thus we deduce 
	\begin{align}\label{eq:lowDeltalarge2}
		& \int_{\theta_0+\frac\delta 2}^{\theta_0+\pi -\frac\delta 2}\int_{\theta_0+\frac\delta 2}^{\theta_0+\pi -\frac\delta 2} \mathbf 1_{|t-s|<\delta} |\alpha(t)-\alpha(s)|\, dtds
		\nonumber
		\\
		&\geq  \delta^2 D\alpha([\theta_0+ 3\delta/2,\theta_0+\pi-3\delta/2]).
	\end{align}
	Similarly we write
	\begin{align}\label{eq:lowDeltalarge3}
		& \int_{\theta_0 -\delta}^{\theta_0 +\delta}
		\int_{\theta_0 -\delta}^{\theta_0 +\delta}
		|\alpha(t)-\alpha(s)|\, dtds \nonumber\\
		&=\int_{\theta_0 -\delta}^{\theta_0 +\delta}
		\left( \iint_{[\theta_0-\delta,\theta_0+\delta]^2} (\mathbf 1_{t<\tau<s} +\mathbf 1_{s<\tau<t})\, dtds \right)\, D\alpha(d\tau), \nonumber\\
		&\leq 
		4\delta^2 D\alpha([\theta_0-\delta,\theta_0+\delta]).
	\end{align}
	Plugging \eqref{eq:lowDeltalarge2}-\eqref{eq:lowDeltalarge3} into \eqref{eq:lowDeltalarge1} we obtain
	\begin{align*}
		\Delta_\varphi(m_1,m_2)&\geq 4\delta^2\left(\frac{1}{\pi} D\alpha([\theta_0+3\delta/2,\theta_0+\pi-3\delta/2]) - 2 D\alpha([\theta_0-\delta,\theta_0+\delta])\right),
	\end{align*}
	if $\pi-\delta <\dist_{\mathbb{S}^1}(m_1,m_2)\leq\pi$.
	Since $\alpha$ is continuous, the measure $D\alpha$ has no atoms, and we deduce the convergence
	\begin{align*}
		&\frac{1}{\pi} D\alpha([\theta_0+3\delta/2,\theta_0+\pi-3\delta/2]) - 2 D\alpha([\theta_0-\delta,\theta_0+\delta])\\
		&\longrightarrow \frac{1}{\pi}D\alpha([\theta_0,\theta_0+\pi])=1,
	\end{align*}
	as $\delta\to 0$, uniformly with respect to $\theta_0\in\R$. In particular we may choose $\delta\in(0,\delta_0)$ sufficiently small such that
	\begin{align*}
		\Delta_\varphi(m_1,m_2)&\geq \delta^2 \geq \frac{\delta^2}{\sup_{\mathbb S^1\times S^1}\Lambda}\Lambda(m_1,m_2),
	\end{align*}
	for $\pi-\delta <\dist_{\mathbb{S}^1}(m_1,m_2)\leq\pi$, and this concludes the proof of Lemma~\ref{l:lowerDelta}.
\end{proof}

\begin{lem}\label{l:lowerLambda}
	For all $m_1,m_2\in\mathbb S^1$, the function $\Lambda$ defined in \eqref{eq:Lambda1} satisfies
	\begin{align*}
		\Lambda(m_1,m_2)\gtrsim \delta^2\omega^{-1}\lt(\delta/2\rt),\qquad\delta=|m_1-m_2|,
	\end{align*}
	where $\omega$ is the minimal modulus of continuity of $\alpha^{-1}$.
\end{lem}
\begin{proof}
	By definition of the modulus of continuity $\omega$ we have
	\begin{align*}
		|\alpha(t)-\alpha(s)|\geq \omega^{-1}(|t-s|)\qquad \forall s,t\in\R.
	\end{align*}
	Using this and the fact that $\omega^{-1}$ is increasing in the definition \eqref{eq:Lambda1} of $\Lambda$ we obtain, for $|\theta_2-\theta_1|\leq\pi$,
	\begin{align*}
		\Lambda(e^{i\theta_1},e^{i\theta_2})& \geq \int_{\theta_1}^{\theta_2}\int_{\theta_1}^{\theta_2} \omega^{-1}(|t-s|)\, dtds\\
		& \geq  \int_{\theta_1}^{\theta_2}\int_{\theta_1}^{\theta_2}\mathbf 1_{\frac 12 |\theta_2-\theta_1| \leq |t-s|\leq |\theta_1-\theta_2|} \, dtds\, \omega^{-1}\left(\frac{|\theta_2-\theta_1|}{2}\right)\\
		&\gtrsim |\theta_2-\theta_1|^2\, \omega^{-1}\left( \frac{|\theta_2-\theta_1|}{2}\right).
	\end{align*}
	Since $\delta=|e^{i\theta_1}-e^{i\theta_2}|\leq|\theta_1-\theta_2|$  and $\omega^{-1}$ is increasing we deduce
	\begin{align*}
		\Lambda(m_1,m_2)&\gtrsim \delta^2 \,\omega^{-1}\lt(\delta/2 \rt),\qquad \delta=|m_1-m_2|,
	\end{align*}
	for all $m_1,m_2\in\mathbb S^1$.
\end{proof}

%
%
%

\subsection{Regularity implies finite entropy production for analytic norms}
\label{ss:revreg}

Recall the definition of $\alpha$ in \eqref{eq:alpha}, i.e. $\alpha\in C^0(\R)$ satisfies
\begin{equation*}
	e^{i\alpha(\theta)} = \gamma'(\theta)\qd\qd\forall \theta\in \R,
\end{equation*}
and that $\alpha$ is increasing and therefore $D\alpha$ is a nonnegative measure on $\R$. 
Here we prove Theorem \ref{thmbesrev}, which follows from the following Lemmas \ref{LBV1} and \ref{LD1}.
\begin{lem} 
	\label{LBV1}
	Let $m$ satisfy \eqref{eq:geneikon}. Assume $D\alpha$ forms a doubling measure and for any $\Omega'\subset \subset \Omega$ we have
	\begin{align}
		\label{eqbe2}
		\sup_{\lt|h\rt|<\mathrm{dist}(\Omega', \partial\Omega)}
		\frac{1}{\lt|h\rt|}\int_{\Omega'} \Pi\lt(m(x),m(x+h)\rt)\, dx <\infty,
	\end{align}
	where $\Pi$ is defined by \eqref{eq:Lambda}. Then the entropy productions of $m$ satisfy $\na\cdot\Phi(m)\in\mathcal{M}_{loc}(\Omega)$ for all $\Phi\in\mathrm{ENT}$, and their lowest upper bound measure  satisfies the estimate
	\begin{align*}
		\left(\bigvee_{\|\lambda_\Phi'\|_\infty\leq 1}
		\lt| \nabla\cdot\Phi(m) \rt|\right)(\Omega')
		& \leq C\,
		\sup_{\lt|h\rt|<\mathrm{dist}(\Omega', \partial\Omega)}
		\frac{1}{\lt|h\rt|}\int_{\Omega'} \Pi\lt(m(x),m(x+h)\rt)\, dx,
	\end{align*}
	for some constant $C>0$ depending on $\Bs$ and the doubling constant of $D\alpha$.
\end{lem}

\begin{proof}   
	Let $m_{\e} = m\ast \rho_{\e}$ for a regularizing kernel $\rho_{\e}$. For any $\Phi\in\mathrm{ENT}$ and any test function $\zeta\in C_c^{\infty}(\Omega)$ with support inside an open subset $V\subset\Omega$, using estimates similar to those leading to \eqref{eq:dvPhimzeta}, we deduce that
	\begin{align}\label{eq:dvPhimzeta2}
		\lt|\langle\nabla\cdot\Phi(m),\zeta\rangle\rt|  \lesssim \|\nabla\Psi\|_\infty \|\zeta\|_\infty \liminf_{\ep\to\infty}\int_V \lt|D m_{\ep}\rt| \lt(1-\|m_{\ep}\|^2\rt)\,dx.
	\end{align}
	Given $x\in V$, note that 
	\begin{align}
		\label{apxeq7}
		\lt|D m_{\ep}(x)\rt|&\leq \int_{B_{\ep}(x)} \lt|m\lt(z\rt)-m\lt(x\rt)\rt|\lt|\na \rho_{\ep}\lt(x-z\rt)\rt| dz\nn\\
		&\lesssim \ep^{-1} \Xint{-}_{B_{\ep}(x)} \lt|m\lt(z\rt)-m\lt(x\rt)\rt| dz.
	\end{align}
	Define the function $F:\R^2\to\R$ by $F(z)=\|z\|^2$ for any $z\in \mathbb{R}^2$. By  convexity of $F$, we have
	\begin{align*}
		1-\|m_{\ep}(x)\|^2&=
		F\lt(m(x)\rt)-F\lt(m_{\ep}(x)\rt) \\
		&\leq 
		\na F\lt(m(x)\rt)\cdot\lt(m(x)-m_{\ep}(x)\rt)
		\\
		&
		=  \na F\lt(m(x)\rt)\cdot \int_{B_{\ep}(x)} \lt(m(x)-m(z)\rt)\rho_{\ep}(x-z)  dz .
	\end{align*}	
	As the level sets $\lbrace F=\lambda^2\rbrace$ are the curves $\lbrace \lambda\gamma(\theta)\rbrace_{\theta\in\R}$, the gradient of $F$ at $m=\gamma(\theta)$ is in the direction of $-i\gamma'(\theta)$. Since moreover $F$ is locally Lipschitz we have 
	\begin{align*}
		\nabla F(\gamma(\theta))=-c(\theta)i\gamma'(\theta),\qquad 0 < c(\theta) \leq C,
	\end{align*}
	for some constant $C>0$ depending on $\Bs$. (Explicitly, $c(\theta)=2/(i\gamma(\theta)\cdot\gamma'(\theta))$.) 
	We write
	$m(x)=\gamma(\theta(x))$ for some $\theta(x)\in\R$ and $m(z)=\gamma(\theta_x(z))$ for some $\theta_x(z)\in\R$ such that $\dist_{\partial \Bs}(m(x),m(z))=|\theta(x)-\theta_x(z)|$, where $\dist_{\partial \Bs}$ denotes the geodesic distance in $\partial\Bs$, and plug the above expression for $\nabla F$ into the previous inequality. This yields
	\begin{align*}
		1-\|m_{\ep}(x)\|^2&\leq C 
		(-i\gamma'(\theta(x)))\cdot 
		\int_{B_{\ep}(x)}  \int_{\theta_x(z)}^{\theta(x)} \gamma'(s)\, ds\,\rho_{\ep}(x-z)  dz\\
		&= C \,
		\int_{B_{\ep}(x)} 
		\int_{\theta_x(z)}^{\theta(x)} 
		\sin(\alpha(\theta(x))-\alpha(s))
		\,ds
		\,\rho_{\ep}(x-z)  dz.
	\end{align*}
	For the last equality we used the definition of the continuous increasing function $\alpha$ characterized by $\gamma'=e^{i\alpha}$. Letting $g(x,z)=\dist_{\partial \Bs}(m(x),m(z))$, we infer
	\begin{align}\label{apxeq4.5}
		1-\|m_{\ep}(x)\|^2&\leq 
		C \,
		\int_{B_{\ep}(x)} 
		\int_{\theta(x)-g(x,z)}^{\theta(x)+g(x,z)} 
		|\alpha(\theta(x))-\alpha(s)|
		\,ds
		\,\rho_{\ep}(x-z)  dz.
	\end{align}
	For all $\theta, r\in\R$, define 
	\begin{align}
		\label{apxeq80.3}
		\GI_{\theta}(r)=\int_{\theta-r}^{\theta+r}\int_{\theta-r}^{\theta+r} \lt|\alpha(t)-\alpha(s)\rt| dt ds.
	\end{align}
	It follows that
	\begin{align}
		\label{apxeq81.3}
		\GI'_{\theta}(r)=2\int_{\theta-r}^{\theta+r}\lt(\lt|\alpha\lt(\theta+r\rt)-\alpha(s)\rt|+ \lt|\alpha\lt(\theta-r\rt)-\alpha(s)\rt|\rt) ds.
	\end{align}
	Using that $\alpha$ is strictly increasing, we find that $\GI'_\theta$ is  strictly increasing, and thus $\GI_\theta$ is strictly convex.

	Moreover, for $\theta-r<s<\theta$ we have $|\alpha(\theta)-\alpha(s)|<|\alpha(\theta+r)-\alpha(s)|$, and for $\theta <s<\theta+r$ we have 
	$|\alpha(\theta)-\alpha(s)|<|\alpha(\theta-r)-\alpha(s)|$. So we deduce from the estimate \eqref{apxeq4.5} and the expression \eqref{apxeq81.3} of $\GI_\theta'$ that
	\begin{align*}
		1-\|m_{\ep}(x)\|^2  &\leq \frac{C}{2}\int_{B_{\ep}(x)}\GI_{\theta(x)}'\lt(g(x,z)\rt)\rho_{\ep}(x-z) dz. 
	\end{align*}
	Putting this together with \eqref{apxeq7}, we obtain 
	\begin{equation*}
		\lt|Dm_{\ep}(x)\rt|\lt(1-\|m_{\ep}(x)\|^2\rt)\lesssim \frac{C}{\e}\,\Xint{-}_{B_{\ep}(x)} \Xint{-}_{B_{\ep}(x)}   \GI_{\theta(x)}'\lt(g(x,z)\rt)  g(x,y)     \, dz\, dy.
	\end{equation*}
	Let $\HI_{\theta}$ denote the Legendre transform of $\GI_{\theta}$, i.e. $\HI_{\theta}(p) = \sup_{r\in\R}\{pr-\GI_{\theta}(r)\}$ for all $p\in\R$. It follows that
	\begin{align}
		\label{eq:product}
		\frac{\e}{C}\,	\lt|Dm_{\ep}(x)\rt|\lt(1-\|m_{\ep}(x)\|\rt)&\lesssim  \Xint{-}_{B_{\ep}(x)} \HI_{\theta(x)}\lt(\GI_{\theta(x)}'\lt(g(x,z)\rt)\rt)dz\nn\\
		&\qd  +  \Xint{-}_{B_{\ep}(x)} \GI_{\theta(x)}\lt(g(x,y)\rt)dy.
	\end{align}
	Note that $\GI_{\theta}(r)\geq cr^2$ for all $r$ sufficiently large and for some universal constant $c>0$. Therefore, for all $p\in\R$, we have $\HI_{\theta}(p) = pr^*-\GI_{\theta}(r^*)$ for the unique $r^*\in\R$ characterized by $\GI_\theta'(r^*)=p$. Thus for all $\theta,r\in\R$ we have
	\begin{align*}
		\HI_{\theta}(\GI_{\theta}'(r))&=\GI_{\theta}'(r)r-\GI_{\theta}(r)\\
		&\leq \GI_{\theta}'(r)r =2r \int_{\theta-r}^{\theta+r}\lt(\lt|\alpha\lt(\theta+r\rt)-\alpha(s)\rt|+ \lt|\alpha\lt(\theta-r\rt)-\alpha(s)\rt|\rt) ds \\
		&\leq 8\, r^2\, |\alpha(\theta+r)-\alpha(\theta-r)|
		=8\, r^2\, D\alpha\lt([\theta-r, \theta+r]\rt).
	\end{align*}
	For the last inequality we used again the fact that $\alpha$ is increasing.
	On the other hand, it is clear from \eqref{apxeq80.3} that
	\begin{equation*}
		\GI_{\theta}\lt(r\rt)\leq 4\, r^2\, D\alpha\lt([\theta-r, \theta+r]\rt).
	\end{equation*}
	Plugging these two inequalities for $\HI_\theta(\GI_\theta')$ and $\GI_\theta$ into \eqref{eq:product} and changing $y$ to $z$ in the second integral on the right-hand side we obtain	
	\begin{align}\label{eq:product2}
		&\lt|Dm_{\ep}(x)\rt|\lt(1-\|m_{\ep}(x)\|\rt)
		\nonumber		
		\\
		&
		\lesssim \frac{C}{\e}\,\Xint{-}_{B_{\ep}(x)}g(x,z)^2D\alpha\lt([\theta(x)-g(x,z), \theta(x)+g(x,z)]\rt)dz. 
	\end{align} 
	Recall from  Remark~\ref{remlam} that, for any $\theta_1,\theta_2\in\R$ with $r=|\theta_1-\theta_2|\leq\pi$, we have
	\begin{align*}
		\Pi(\gamma(\theta_1),\gamma(\theta_2))&=
		\Lambda(e^{i\theta_1},e^{i\theta_2})
		\gtrsim r^2 D\alpha\left(\left[\frac{\theta_1+\theta_2}{2}-\frac r6, 
		\frac{\theta_1+\theta_2}{2}+\frac r6 \right]\right).
	\end{align*}
	Using the fact that $D\alpha$ is a doubling measure, we deduce
	\begin{align*}
		\Pi(\gamma(\theta_1),\gamma(\theta_2))&\gtrsim C_0\, r^2 D\alpha\left(\left[\frac{\theta_1+\theta_2}{2}-\frac{3 r}{2}, 
		\frac{\theta_1+\theta_2}{2}+\frac{3 r}{2} \right]\right)\\
		&\geq C_0 \, r^2 D\alpha\lt([\theta_1-r, \theta_1+r]\rt),
	\end{align*}	
	for some constant $C_0>0$ depending on the doubling constant of $D\alpha$. Applying this to $\theta_1=\theta(x)$, $\theta_2=\theta_x(z)$ and $r=g(x,z)$ and plugging the resulting inequality into \eqref{eq:product2} we find
	\begin{align*}
		\lt|Dm_{\ep}(x)\rt|\lt(1-\|m_{\ep}(x)\|\rt)
		\leq \frac{C_1}{\e}\,\Xint{-}_{B_\e(x)}\Pi(m(x),m(z)) \, dz,
	\end{align*}
	for a constant $C_1>0$ depending on the doubling constant of $D\alpha$  and $\Bs$.
	
	Integrating this estimate with respect to $x\in V$ and recalling \eqref{eq:dvPhimzeta2} we deduce
	\begin{align*}
		\lt|\langle\nabla\cdot\Phi(m),\zeta\rangle\rt| 
		\lesssim C_1 \|\nabla\Psi\|_\infty \|\zeta\|_\infty 
		\liminf_{\e\to 0}
		\frac{1}{\e}\Xint{-}_{B_\e(0)}
		\int_V \Pi(m(x),m(x+h))\, dx\, dh,
	\end{align*} 
	for any $\zeta\in C_c^\infty(V)$. Thanks to \eqref{eqbe2}, the limit in the right-hand side is finite. This implies in particular that $\nabla\cdot\Phi(m)$ is a locally finite Radon measure, such that
	\begin{align*}
		\lt|\nabla\cdot\Phi(m) \rt|(V)\lesssim C_1 
		\|\nabla\Psi\|_\infty  
		\liminf_{\e\to 0}
		\frac{1}{\e}\Xint{-}_{B_\e(0)}
		\int_V \Pi(m(x),m(x+h))\, dx\, dh.
	\end{align*}
	Moreover, from this estimate we infer (arguing as in the proof of Proposition~\ref{p:lowerbound}) that 
	\begin{align*}
		\left(\bigvee_{\|\lambda_\Phi'\|_\infty\leq 1}
		\lt| \nabla\cdot\Phi(m) \rt|\right)(\Omega')
		&  \lesssim  C_1\, \liminf_{\e\to 0}
		\frac{1}{\e}\Xint{-}_{B_\e(0)}
		\int_{\Omega'} \Pi(m(x),m(x+h))\, dx\, dh,
	\end{align*}
	for any open subset $\Omega'\subset\subset\Omega$. This implies the conclusion of Lemma~\ref{LBV1}.
\end{proof}
%
%
%

\begin{lem}
\label{LD1}
Suppose $\partial \Bs$ is analytic. Then $\alpha$ is analytic and $\alpha'(t)\,dt$ forms a doubling measure.
\end{lem}

\begin{proof} The lemma comes down to the fact that an absolutely continuous measure whose density is a nonnegative analytic function is doubling. This is presumably a well known fact, but we found no direct reference for it beyond a (more general) theorem in \cite{gar} for the square of an analytic function. Since $\alpha'\geq 0$, the function $\beta(t):=\sqrt{\alpha'(t)}$ is well defined. As the square root function is analytic in $(0,\infty)$, it follows that $\beta(t)$ is analytic at all $t$ such that $\alpha'(t)>0$ (see \cite[Proposition 1.4.2]{krantz}). Given $t_0$ with $\alpha'(t_0)=0$, again by the fact that $\alpha'\geq 0$, we can write, for $t$ in a sufficiently small neighborhood $I$ of $t_0$, $\alpha'(t)=(t-t_0)^{2p} h(t)$ for some integer $p\geq 1$ and some analytic function $h$ with $h(t)\ne 0$ in $I$. Thus $\beta(t)=(t-t_0)^p\sqrt{h(t)}$ is analytic in $I$. This shows that $\beta(t)$ is analytic in $\R$, and thus \cite[Theorem 1]{gar} can be applied to $\beta$ to conclude that $\alpha'(t)\,dt=\beta^2(t)\,dt$ is a doubling measure.
\end{proof}

\begin{rem}
	\label{rem:lp-sphere}
	The unit sphere of the $\ell^p$ norm in $\R^2$ defined by $|x|^p+|y|^p=1$ is analytic for $1<p<\infty$. This follows from the analyticity of the function $\lt(1-|x|^p\rt)^{1/p}$ in $(-1,1)$, which in turn is a consequence of the analyticity of the functions $1-|x|^p$ and $x^{1/p}$ in the intervals $(-1,1)$ and $(0,\infty)$, respectively.
\end{rem}

%
%

\section{Comparison of upper and lower bounds when $m$ is $BV$}\label{s:boundsBV}

%
%

In this section we consider $m\in BV(\Omega;\R^2)$ that satisfies \eqref{eq:geneikon} and assume that $\Omega\subset\R^2$ is a bounded simply connected smooth domain. Therefore the constraint $\nabla\cdot m=0$ is equivalent to the existence of a function $u$ such that $m=i\nabla u$. Using that correspondence, it is  somewhat lengthy but straightforward to see that the upper bound construction in \cite{ark2} directly applies (taking $F(A,w)=\lt|A\rt|^2+\lt(1-\|iw\|^2\rt)^{2}$ for $A\in\R^{2\times 2}$ and $w\in\R^2$ in  \cite[Theorem~1.2]{ark2}) to provide the existence of a $C^1$ sequence $m_\e\to m$ in $L^p(\Omega;\R^2)$ for $1\leq p <\infty$, such that $\nabla\cdot m_{\e}=0$ and
\begin{align*}
\limsup_{\e\to 0} I_\e(m_\e) \leq \int_{J_m} \mathrm{c}^{\mathrm{1D}}(m^+,m^-)\, d\mathcal H^1.
\end{align*}
Here $m^\pm$ are the traces of $m$ along its jump set $J_m$, and $\mathrm{c}^{\mathrm{1D}}(z^+,z^-)$ is the optimal energy of a one-dimensional transition between two states $z^\pm\in\partial\Bs$. In other words
\begin{align}\label{eq:c1Dinf}
\mathrm{c}^{\mathrm{1D}}(z^+,z^-)&=\inf_{\zeta\in Y} \left\lbrace \int_{-\infty}^{+\infty} \left( (\zeta'(x))^2 + (1-\|a\,\nu +\zeta(x) \, i\nu\|^2)^2\right)\, dx\right\rbrace,
\\
\text{where }
Y&=\left\lbrace \zeta\in C^1(\R)\colon \lim_{x\to\pm\infty}\zeta =z^\pm\cdot i\nu\right\rbrace,
\nonumber
\end{align}
and
\begin{align}\label{eq:nuzpm}
\nu &=\nu_{z^+,z^-}=i\frac{z^+-z^-}{|z^+-z^-|},\qquad a=a_{z^+,z^-}=z^+\cdot\nu=z^-\cdot\nu.
\end{align}
Classically (see e.g. \cite{sternberg91} for details), the infimum in \eqref{eq:c1Dinf} can be explicitly calculated. Indeed, assuming without loss of generality that $z^-\cdot i\nu \leq z^+\cdot i\nu$, for any admissible function $\zeta\in Y$ we have
\begin{align*}
&\int_{-\infty}^{+\infty} \left( (\zeta'(x))^2 + (1-\|a\,\nu +\zeta(x) \, i\nu\|^2)^2\right)\, dx \\
&
\geq 2\int_{-\infty}^{+\infty} (1-\|a\,\nu +\zeta(x)\,i\nu\|^2) \, \zeta'(x) \, dx \\
&= 2 \int_{-\infty}^{+\infty} \frac{d}{dx}\left[\int_{-\infty}^{\zeta(x)} (1-\|a\,\nu + s \, i\nu\|^2)\, ds\right]\, dx\\
&=2\int_{z^-\cdot i\nu}^{z^+\cdot i\nu} (1-\|a\,\nu + s \, i\nu\|^2)\, ds,
\end{align*}
and conversely, one can check that any solution of $\zeta'=1-\|a\,\nu +\zeta\, i\nu\|^2$ with initial condition $\zeta(0)\in (z^-\cdot i\nu,z^+\cdot i\nu)$ is admissible, i.e. belongs to the class $Y$, and achieves equality. So we have
\begin{align*}
\mathrm{c}^{\mathrm{1D}}(z^+,z^-)=
2\left|\int_{z^{-}\cdot i\nu}^{z^{+}\cdot i\nu} \lt(1-\|a \nu + s\, i \nu \|^2\rt) ds\right|,
\end{align*}
which corresponds to the expression \eqref{eq:c1D} given in the introduction.

Our goal in this section is to prove Theorem \ref{t:lowerupper} by comparing this upper bound with the lower bound \eqref{eq:lowerbound} provided by the entropy productions:
\begin{align*}
\left(\bigvee_{\|\lambda_\Phi'\|_{\infty}\leq 1}|\nabla\cdot\Phi(m)|\right)(\Omega) \leq C_0 \,\liminf_{\e\to 0} I_{\e}(m_{\e}).
\end{align*}
This follows from the estimate \eqref{eq:cENTleqc1D} and Lemma \ref{l:c1DleqcENT} below.

\begin{lem}\label{l:cENT}
For $m\in BV(\Omega;\R^2)$ satisfying \eqref{eq:geneikon} we have
\begin{align*}
\bigvee_{\|\lambda_\Phi'\|_{\infty}\leq 1}|\nabla\cdot\Phi(m)| = \mathrm{c}^{\mathrm{ENT}}(m^+,m^-)\, \mathcal H^1_{\lfloor J_m},
\end{align*}
where
\begin{align}\label{eq:cENT}
 \mathrm{c}^{\mathrm{ENT}}(z^+,z^-)=\sup_{\lambda\in\Lambda_*} \left\lbrace \int_{\theta^-}^{\theta^+}\lambda(t)\,\gamma'(t)\cdot\nu_{z^+,z^-}\, dt\right\rbrace \qquad\text{for }z^\pm=\gamma(\theta^\pm),
\end{align}
and
\begin{align*}
\Lambda_*=\left\lbrace\lambda\in C^1(\R/2\pi\Z)\colon \int_{\R/2\pi \Z}\lambda(t)\gamma'(t)\, dt = 0\text{ and }\|\lambda'\|_\infty\leq 1\right\rbrace.
\end{align*}
\end{lem}
\begin{proof}
For any $\Phi\in \mathrm{ENT}$, the function $\lambda_\Phi$ defined by $\frac{d}{d\theta}\Phi(\gamma(\theta))=\lambda_\Phi(\theta)\gamma'(\theta)$ belongs to $\Lambda_*$ if it satisfies $\|\lambda_\Phi'\|_\infty\leq 1$, and reciprocally, to any $\lambda\in\Lambda_*$ 
one can associate an entropy
 $\Phi_\lambda\in \mathrm{ENT}$ by setting
\begin{align*}
\Phi_\lambda(\gamma(\theta))=\int_0^\theta \lambda(t)\, \gamma'(t)\, dt.
\end{align*}
With these notations we therefore have
\begin{align*}
\bigvee_{\|\lambda_\Phi'\|_{\infty}\leq 1}|\nabla\cdot\Phi(m)|
=\bigvee_{\lambda\in\Lambda_*}|\nabla\cdot\Phi_{\lambda}(m)|.
\end{align*}
For a $BV$ map $m$,  the $BV$ chain rule implies that the entropy productions are absolutely continuous with respect to $\mathcal H^1_{\lfloor J_m}$, and
\begin{align*}
|\dv\Phi_\lambda(m)|&= |(\Phi_{\lambda}(m^+)-\Phi_{\lambda}(m^-))\cdot\nu | \, d\mathcal H^1_{\lfloor J_m}\\
&=|\mathrm{c}_{\lambda}(m^+,m^-)|\, \, d\mathcal H^1_{\lfloor J_m},\\
\mathrm{c}_{\lambda}(z^+,z^-)&=\int_{\theta^-}^{\theta^+}\lambda(t)\,\gamma'(t)\cdot \nu_{z^+,z^-}\, dt\quad \text{for }z^\pm=\gamma(\theta^\pm).
\end{align*}
Therefore, restricting the supremum to a countable dense subset of $\Lambda_*$ and applying \cite[Remark~1.69]{ambrosio}, we see that the lowest upper bound measure is also absolutely continuous with respect to $\mathcal H^1_{\lfloor J_m}$ and given by
\begin{align*}
\bigvee_{\lambda\in\Lambda_*}|\nabla\cdot\Phi_\lambda(m)|
&=\left(\sup_{\lambda\in\Lambda_*} |\mathrm{c}_\lambda(m^+,m^-)|\right)\, \mathcal H^1_{\lfloor J_m}.
\end{align*}
Since $\lambda\mapsto\mathrm{c}_\lambda$ is linear we can remove the absolute value in the right-hand side, concluding the proof of Lemma~\ref{l:cENT}.
\end{proof}

Combining the lower and upper bounds, we see that for any $m\in BV(\Omega;\R^2)$ satisfying \eqref{eq:geneikon} we have
\begin{align*}
\int_{J_m}\mathrm{c^{\mathrm{ENT}}}(m^+,m^-)\,d\mathcal H^1 \leq C_0\int_{J_m}\mathrm{c}^{\mathrm{1D}}(m^+,m^-)\, d \mathcal H^1.
\end{align*}
For any fixed $z^\pm\in\partial\Bs$ we may apply this to a divergence-free map taking only the two values $z^\pm$, and we immediately deduce the inequality
\begin{align}\label{eq:cENTleqc1D}
\mathrm{c}^{\mathrm{ENT}}(z^+,z^-)\leq C_0\, \mathrm{c}^{\mathrm{1D}}(z^+,z^-)\qquad\forall z^\pm\in\partial\Bs.
\end{align}
Next we prove the reverse inequality. To that end we start by obtaining a more explicit expression of $\mathrm{c}^{\mathrm{ENT}}$ for small jumps.

\begin{lem}
	\label{LL11.3}
	For all $z^{\pm}=\gamma(\theta^\pm) \in \partial\Bs$ with $|\theta^+-\theta^-|=\mathrm{dist}_{\partial\Bs}(z^+,z^-)<\pi/2$, we have 
	\begin{equation}
		\label{eqaqq82}
		\mathrm{c}^{\mathrm{ENT}}\lt(z^{+},z^{-}\rt) =\lt|\int_{\theta^-}^{\theta^+}\lt(\theta-\ti\theta_{z^+,z^-}\rt)\lt(\gamma'(\theta)\cdot \nu_{z^+,z^-}\rt)\,d\theta\rt|,
	\end{equation}
	where 
	$\ti\theta_{z^+,z^-}\in\R$ is the unique point between $\theta^-$ and $\theta^+$ satisfying
\begin{equation*}
	\gamma'\lt(\ti\theta_{z^+,z^-}\rt)\cdot \nu_{z^+,z^-}=0.
\end{equation*}
Moreover we also have
	\begin{equation}
		\label{eqkk11.22}
		\inf_{\mathrm{dist}_{\partial\Bs}(z^+,z^-)\geq\frac\pi 2}\mathrm{c}^{\mathrm{ENT}}\lt(z^{+},z^{-}\rt) >0.
	\end{equation}
\end{lem}
\begin{proof} 
Let $z^{\pm}=\gamma(\theta^\pm) \in \partial \Bs$ with $|\theta^+-\theta^-|<\pi/2$. Assume without loss of generality that $\theta^- < \theta^+  < \theta^-+\pi/2$.
To simplify notations, in this proof we drop the indices $(z^+,z^-)$ and simply write $\nu=\nu_{z^+,z^-}$ and $\tilde\theta=\tilde\theta_{z^+,z^-}$.

Recall that $\mathrm{c}^{\mathrm{ENT}}$ is given by
\begin{align*}
\mathrm{c}^{\mathrm{ENT}}(z^+,z^-)&=\sup_{\lambda\in\Lambda_*} \mathrm{c}_\lambda(z^+,z^-),\\
\text{where }\mathrm{c}_\lambda(z^+,z^-)&=\int_{\theta^-}^{\theta^+}\lambda(\theta)\gamma'(\theta)\cdot \nu \,d\theta\\
&=\int_{\theta^-}^{\theta^+}\lt(\lambda(\theta)-\lambda(\ti\theta )\rt)\gamma'(\theta)\cdot \nu \,d\theta.
\end{align*}
The last equality is valid because $\int_{\theta^-}^{\theta^+}\gamma'(\theta)\cdot\nu\, d\theta=(z^+-z^-)\cdot\nu=0$ by definition of $\nu$. Since  for 
all $\lambda\in\Lambda_*$ we have $\|\lambda'\|_{\infty}\leq 1$, and therefore 
$|\lambda (\theta)-\lambda (\tilde\theta)|\leq |\theta-\tilde\theta|$. This implies
\begin{align*}
\mathrm{c}^{\mathrm{ENT}}(z^+,z^-)&\leq \int_{\theta^-}^{\theta^+}
|\theta-\ti\theta |\,\lt|\gamma'(\theta)\cdot \nu \rt|\,d\theta.
\end{align*}
The definition of $\tilde\theta$ and convexity of $\Bs$  ensure that $(\theta-\tilde\theta)(\gamma'(\theta)\cdot\nu)\geq 0$   for $\theta\in (\theta^-,\theta^+)$, so the above becomes
\begin{align*}
\mathrm{c}^{\mathrm{ENT}}(z^+,z^-)&\leq \lt|\int_{\theta^-}^{\theta^+}
(\theta-\ti\theta )\,(\gamma'(\theta)\cdot \nu) \,d\theta\rt|.
\end{align*}
Conversely, since $|\theta^+-\theta^-|<\pi/2$ we may choose a $\pi$-periodic $\lambda_0\in C^1(\R/\pi\Z)$ with $\|\lambda_0'\|_{\infty}\leq 1$ and such that $\lambda_0(\theta)-\lambda_0(\tilde\theta)=\theta-\tilde\theta$ for $\theta^-<\theta<\theta^+$. Note that the $\pi$-periodicity of $\lambda_0$ ensures $\int_{\R/2\pi\Z}\lambda_0\gamma'=0$ since $\gamma'(t+\pi)=-\gamma'(t)$, so $\lambda_0\in\Lambda_*$. Therefore we have $\mathrm{c}^{\mathrm{ENT}}\geq\mathrm{c}_{\lambda_0}$ and
 we deduce that  $\mathrm{c}^{\mathrm{ENT}}$ is given by \eqref{eqaqq82}.

To prove \eqref{eqkk11.22}, note that $\mathrm{c}^{\mathrm{ENT}}$ is defined in \eqref{eq:cENT} as a supremum of continuous functions, and is therefore lower semicontinuous  on $\partial\Bs\times\partial\Bs$. Hence the infimum in \eqref{eqkk11.22} is attained at some $z^\pm\in \partial \Bs$, $\mathrm{dist}_{\partial\Bs}\lt(z^{+},z^{-}\rt)\geq \frac{\pi}{2}$.
 As $\Bs$ is strictly convex, the function $\theta\mapsto \gamma'(\theta)\cdot\nu$ cannot be identically zero on any open interval, which prevents $\mathrm{c}^{\mathrm{ENT}}$ from vanishing unless $z^+=z^-$. So the infimum in \eqref{eqkk11.22} is positive.
\end{proof}

\begin{lem}\label{l:c1DleqcENT}
We have
\begin{align*}
\mathrm{c}^{\mathrm{1D}}(z^+,z^-)\leq C\, \mathrm{c}^{\mathrm{ENT}}(z^+,z^-)\qquad\forall z^\pm\in\partial \Bs,
\end{align*}
for some constant $C>0$ depending only on $\Bs$.
\end{lem}
\begin{proof}
Let $z^\pm=\gamma(\theta^\pm)\in\partial\Bs$ be two distinct points, with $|\theta^+-\theta^-|=\dist_{\partial \Bs}(z^+,z^-)$, and $\theta^-<\theta^+\leq \theta^- +\pi$. Dropping the indices $(z^+,z^-)$ we denote
\begin{align*}
\nu=i\frac{z^+-z^-}{|z^+-z^-|},\qquad a=z^+\cdot\nu=z^-\cdot\nu\leq 0,
\end{align*}
and recall that $\mathrm{c}^{\mathrm{1D}}$ is given by
\begin{align*}
\mathrm{c}^{\mathrm{1D}}(z^+,z^-)=2\int_{z^+\cdot i\nu}^{z^-\cdot i\nu}\left(1-\|a\,\nu + s\,i\nu\|^2\right)\, ds.
\end{align*}
Since $\Bs$ is strictly convex, for any $\theta\in (\theta^-,\theta^+)$ there is  a unique $s(\theta)\in (z^+\cdot i\nu, z^-\cdot i\nu)$ such that
\begin{align*}
a\,\nu  + s(\theta)\,i\nu  = \beta(\theta) \gamma(\theta)\text{ for some }0<\beta(\theta)<1.
\end{align*}
The function $s\colon (\theta^-,\theta^+)\to (z^+\cdot i\nu, z^-\cdot i\nu)$ is a decreasing bijection. Taking the scalar product of the above with $i\gamma(\theta)$ and with $\nu$ we have
\begin{align*}
s(\theta)=a\frac{\gamma(\theta)\cdot i\nu}{\gamma(\theta)\cdot\nu},\qquad \beta(\theta)=\frac{a}{\gamma(\theta)\cdot\nu}.
\end{align*}
The change of variable $s=s(\theta)$ therefore gives
\begin{align}
\mathrm{c}^{\mathrm{1D}}(z^+,z^-)&=-2\int_{\theta^- }^{\theta^+ }\left(1-\beta(\theta)^2\right)\, s'(\theta)\, d\theta
\nonumber\\
&=-2\int_{\theta^- }^{\theta^+ }\left(1-\frac{a^2}{(\gamma(\theta)\cdot\nu)^2}\right)\, s'(\theta)\, d\theta 
\nonumber \\
&=-2\int_{\theta^- }^{\theta^+ }(\gamma(\theta)\cdot\nu -a) \frac{\gamma(\theta)\cdot\nu +a}{(\gamma(\theta)\cdot\nu)^2} \, s'(\theta)\, d\theta.
\label{eq:Achangevar}
\end{align}
Since $|\gamma'|=1$, from the explicit expression of $s(\theta)$ we have $|s'(\theta)|\leq 2 |a| \,|\gamma(\theta)|/(\gamma(\theta)\cdot\nu)^2$. For $|\theta^+-\theta^-| < \pi$, using moreover the inequality $\gamma(\theta)\cdot\nu \leq a< 0$, which follows from the convexity of $\Bs$ and implies in particular $(\gamma(\theta)\cdot\nu)^2\geq a^2$,  we deduce
\begin{align}\label{eq:Aleqinta}
\mathrm{c}^{\mathrm{1D}}(z^+,z^-)& 
\leq \frac{4}{|a|^3}\int_{\theta^- }^{\theta^+ } (a - \gamma(\theta)\cdot\nu)\,
|\gamma(\theta)|\, |\gamma(\theta)\cdot\nu +a|\, d\theta \nonumber\\
&\lesssim \frac{1}{|a|^3}\int_{\theta^- }^{\theta^+ } (a - \gamma(\theta)\cdot\nu)\, d\theta.
\end{align}
The last inequality follows from $\mathrm{diam}(\Bs)\lesssim 1$.
Recalling from Lemma \ref{LL11.3} the definition of $\tilde\theta$ as the unique $\tilde\theta\in (\theta^-,\theta^+)$ such that $\gamma'(\tilde\theta)\cdot\nu=0$, and that $a=z^+\cdot\nu=z^-\cdot\nu$, we write
\begin{align*}
a-\gamma(\theta)\cdot\nu 
&=
\mathbf 1_{\tilde \theta <\theta<\theta^+} \int_{\theta}^{\theta^+}\gamma'(t)\cdot\nu\, dt
-
\mathbf 1_{\theta^-<\theta<\tilde\theta} \int_{\theta^-}^\theta \gamma'(t)\cdot\nu\, dt,
\end{align*}
and
\begin{align}
\int_{\theta^- }^{\theta^+ } (a - \gamma(\theta)\cdot\nu)\, d\theta 
&=\int_{\tilde\theta}^{\theta^+} \int_{\theta}^{\theta^+}\gamma'(t)\cdot\nu\, dt\,d\theta
-\int_{\theta^-}^{\tilde\theta} \int_{\theta^-}^{\theta}\gamma'(t)\cdot\nu\, dt\,d\theta
\nonumber
\\
&=\int_{\tilde\theta}^{\theta^+} \int_{\tilde\theta}^{t}\, d\theta\, \gamma'(t)\cdot\nu\, dt 
-\int_{\theta^-}^{\tilde\theta} \int_{t}^{\tilde \theta}\, d\theta\, \gamma'(t)\cdot\nu\, dt
\nonumber\\
&=\int_{\theta^-}^{\theta^+}(t-\tilde\theta)(\gamma'(t)\cdot\nu)\,dt.
\label{eq:intequalsB}
\end{align}
For $|\theta^+-\theta^-| < \pi/2$, this last integral is exactly the expression of $\mathrm{c}^{\mathrm{ENT}}(z^+,z^-)$ given by Lemma~\ref{LL11.3}. Therefore, combining this with \eqref{eq:Aleqinta} we deduce
\begin{align*}
\mathrm{c}^{\mathrm{1D}}(z^+,z^-)& 
\lesssim \frac{1}{|a_{z^+,z^-}|^3} \mathrm{c}^{\mathrm{ENT}}(z^+,z^-)\quad\text{for }\dist_{\partial\Bs}(z^+,z^-) < \frac\pi 2.
\end{align*}
The function
\begin{align*}
(z^+,z^-)\mapsto |a_{z^+,z^-}|,
\end{align*}
is continuous on $\partial \Bs\times\partial\Bs \setminus \lbrace z^+=z^-\rbrace$, vanishes exactly when $\mathrm{dist}_{\partial\Bs}(z^+,z^-)=\pi$, and for $z^+$ close to $z^-$ it satisfies
\begin{align*}
|a_{z^+,z^-}|\longrightarrow | i\gamma'(\theta)\cdot \gamma(\theta)|,
\quad\text{ as }(z^+,z^-)\to (\gamma(\theta),\gamma(\theta)).
\end{align*} 
By convexity of $\Bs$ we have $| i\gamma'(\theta)\cdot \gamma(\theta)|\geq\alpha_0>0$, where $\alpha_0$ is the largest radius of a euclidean ball contained in $\Bs$. From these properties we deduce that
\begin{align*}
\inf_{\dist_{\partial\Bs}(z^+,z^-) \leq \frac\pi 2 } |a_{z^+,z^-}| >0,
\end{align*}
and the above bound on $\mathrm{c}^{\mathrm{1D}}$ implies
\begin{align*}
\mathrm{c}^{\mathrm{1D}}(z^+,z^-)& 
\leq C \mathrm{c}^{\mathrm{ENT}}(z^+,z^-)\quad\text{for }\dist_{\partial\Bs}(z^+,z^-) < \frac\pi 2,
\end{align*}
for some constant $C>0$ depending only on $\Bs$. Since $\mathrm{c}^{\mathrm{1D}}\leq 2 \pi$ and $\mathrm{c}^{\mathrm{ENT}}$ is bounded away from zero for $\dist_{\partial \Bs}(z^+,z^-)\geq\pi/2$ thanks to Lemma~\ref{LL11.3}, this inequality is true for all $z^\pm\in\partial \Bs$.
\end{proof}

\begin{rem}
Combining the expressions \eqref{eq:Achangevar} and \eqref{eq:intequalsB} obtained in the proof of Lemma~\ref{l:c1DleqcENT} and passing to the limit as $\theta^+-\theta^-\to 0$ one obtains
\begin{align*}
\frac{\mathrm{c}^{\mathrm{1D}}(z^+,z^-)}{\mathrm{c}^{\mathrm{ENT}}(z^+,z^-)}\longrightarrow \frac{2}{|i\gamma'(\theta)\cdot\gamma(\theta)|}\qquad\text{as }(z^+,z^-)\to (\gamma(\theta),\gamma(\theta)).
\end{align*}
Hence for infinitesimally small jumps, the costs $\mathrm{c}^{\mathrm{1D}}$ and $\mathrm{c}^{\mathrm{ENT}}$ differ by the above multiplicative factor, which depends on the direction of the jump.
\end{rem}

\subsection{Regularity controls entropy productions when $m$ is $BV$}
\label{besbv}

In this subsection we prove Theorem~\ref{t:regconvBV}. To that end we first compare the jump cost $c^{\mathrm{ENT}}$ to the regularity cost $\Pi$  defined in \eqref{eq:Lambda}.

\begin{lem}\label{l:cENTleqLambda}
	We have
	\begin{align*}
		\mathrm{c}^{\mathrm{ENT}}(z^+,z^-)\leq  \, \Pi(z^+,z^-),
	\end{align*}
	for all $z^\pm\in \partial\Bs$.
\end{lem}
\begin{proof}
	We let $z^{\pm}=\gamma(\theta^\pm)$ for some $\theta^\pm\in\R$ such that $|\theta^+-\theta^-|=\dist_{\partial\Bs}(z^+,z^-)\leq \pi$, and assume without loss of generality that $\theta^-<\theta^+$.

	From the proof of Lemma~\ref{LL11.3} we have
	\begin{align*}
		\mathrm{c}^{\mathrm{ENT}}(z^+,z^-)
		&\leq \int_{\theta^-}^{\theta^+}(\theta-\tilde\theta)\gamma'(\theta)\cdot\nu\, d\theta,
	\end{align*}
	where $\tilde\theta\in (\theta^-,\theta^+)$ is such that $\nu=i\gamma'(\tilde\theta)$ and $\nu=\nu_{z^+,z^-}$ is defined in \eqref{eq:nuzpm}. 
	Recalling that the continuous increasing function $\alpha$ is defined by $\gamma'=e^{i\alpha}$, we rewrite this as
	\begin{align*}
		\mathrm{c}^{\mathrm{ENT}}(z^+,z^-)
		&\leq \int_{\theta^-}^{\theta^+}(\theta-\tilde\theta)\,e^{i\alpha(\theta)}\cdot ie^{i\alpha(\tilde\theta)}\, d\theta\\
		&=\int_{\theta^-}^{\tilde\theta}(\tilde\theta-\theta)\,\sin(\alpha(\tilde\theta)-\alpha(\theta))\, d\theta\\
		&\quad
		+ 
		\int_{\tilde\theta}^{\theta^+}(\theta-\tilde\theta)\,
		\sin(\alpha(\theta)-\alpha(\tilde\theta))\, d\theta \\
		&\leq (\tilde\theta-\theta^-)A^- +(\theta^+-\tilde\theta) A^+,\\
		\text{where }
		A^-&=\int_{\theta^-}^{\tilde\theta} \sin(\alpha(\tilde\theta)-\alpha(\theta))\, d\theta,
		\quad
		A^+=\int_{\tilde\theta}^{\theta^+} 
		\sin(\alpha(\theta)-\alpha(\tilde\theta))\, d\theta.
	\end{align*}
	Next recall that $(z^+-z^-)\cdot\nu=0$ by definition of $\nu$, and rewrite this as
	\begin{align*}
		0&=\int_{\theta^-}^{\theta^+}\gamma'(\theta)\cdot\nu\, d\theta
		=\int_{\theta^-}^{\theta^+}\sin(\alpha(\theta)-\alpha(\tilde\theta))\, d\theta =A^+-A^-,
	\end{align*}
	so we have in fact $A^+=A^-$ and the above estimate for $\mathrm{c}^{\mathrm{ENT}}$ becomes
	\begin{align}\label{eq:cENTA}
		\mathrm{c}^{\mathrm{ENT}}(z^+,z^-)&\leq (\theta^+-\theta^-) A,
		\nonumber\\
		\text{where } A&=\int_{\theta^-}^{\tilde\theta} \sin(\alpha(\tilde\theta)-\alpha(\theta))\, d\theta=\int_{\tilde\theta}^{\theta^+} 
		\sin(\alpha(\theta)-\alpha(\tilde\theta))\, d\theta.
	\end{align}
	Note that, using the fact that $\alpha$ is absolutely continuous, we have
	\begin{align*}
		A&\leq \int_{\theta^-}^{\tilde\theta}(\alpha(\tilde\theta)-\alpha(\theta))\, d\theta 
		=
		\int_{\theta^-}^{\tilde\theta}\int_{\theta}^{\tilde\theta} D\alpha(d\tau)\, d\theta \
		=\int_{\theta^-}^{\tilde \theta } (\tau-\theta^-)\, D\alpha(d\tau),
	\end{align*}
	and similarly
	\begin{align*}
		A&\leq \int_{\tilde\theta}^{\theta^+}(\theta^+-\tau)\, D\alpha(d\tau).
	\end{align*}
	So from \eqref{eq:cENTA} we infer
	\begin{align}\label{eq:cENTleqmin}
		&\mathrm{c}^{\mathrm{ENT}}(z^+,z^-)
		\nonumber\\
		&\leq (\theta^+-\theta^-) \,
		\min\left \lbrace 
		\int_{\theta^-}^{\tilde \theta } (\tau-\theta^-)\, D\alpha(d\tau),
		\int_{\tilde\theta}^{\theta^+}(\theta^+-\tau)\, D\alpha(d\tau)
		\right\rbrace.
	\end{align}
	
	Next we consider two cases, depending on whether $\tilde\theta\in (\theta^-,\theta^+)$ is closer to $\theta^-$ or to $\theta^+$. If $\tilde\theta$ is closer to $\theta^-$ we have
	$\theta^+-\theta^- \leq 2 (\theta^+-\tilde\theta)$, 
	so from  \eqref{eq:cENTleqmin} (recalling that $\alpha$ is increasing and thus $D\alpha$ is a nonnegative measure) we deduce
	\begin{align*}
		\mathrm{c}^{\mathrm{ENT}}(z^+,z^-)&\leq 2 (\theta^+-\tilde\theta)
		\int_{\theta^-}^{\tilde \theta } (\tau-\theta^-)\, D\alpha(d\tau) \\
		&\leq 2 \int_{\theta^-}^{\tilde \theta }  (\theta^+-\tau) (\tau-\theta^-)\, D\alpha(d\tau).
	\end{align*}
	Otherwise, $\tilde\theta$ is closer to $\theta^{+}$ so we have $\theta^+-\theta^- \leq 2 (\tilde\theta-\theta^-)$ and we find
	\begin{align*}
		\mathrm{c}^{\mathrm{ENT}}(z^+,z^-)
		&\leq 2 
		\int_{\tilde \theta }^{\theta^+} (\theta^+-\tau) (\tau-\theta^-)\, D\alpha(d\tau).
	\end{align*}
	In both cases, we have
	\begin{align*}
		\mathrm{c}^{\mathrm{ENT}}(z^+,z^-)
		&\leq 2 
		\int_{\theta^- }^{\theta^+} (\theta^+-\tau) (\tau-\theta^-)\, D\alpha(d\tau),
	\end{align*}
	and thanks to Remark~\ref{remlam} and the symmetry of $\Pi$, this last expression is exactly $\Lambda(e^{i\theta^-}, e^{i\theta^+})=\Pi(\gamma(\theta^-), \gamma(\theta^+))=\Pi(z^+,z^-)$.
\end{proof}

Next we deduce from Lemma~\ref{l:cENTleqLambda} and properties of $BV$ maps 
that the regularity estimate provided by $\Pi$ controls the entropy productions, proving Theorem~\ref{t:regconvBV}.
Recall for an entropy $\Phi\in \mathrm{ENT}$ the $C^1$ function $\lambda_{\Phi}$ is defined by $\frac{d}{d\theta}\Phi(\gamma(\theta))=\lambda_\Phi(\theta)\gamma'(\theta)$.

\begin{lem}
	\label{misbv}
	Let $m\in BV\lt(\Omega;\R^2\rt)$ satisfy  \eqref{eq:geneikon}. Then for any open set  $\Omega'\subset \subset \Omega$
	we have  
	\begin{align*}
		&\left(\bigvee_{\|\lambda_{\Phi}'\|_{\infty}\leq 1} \lt|\na \cdot \Phi(m)\rt|\right)\lt(\Omega'\rt) 
		\leq C_0 \,\limsup_{\lt|h\rt|\to 0} \frac{1}{|h|}\int_{\Omega'} 
		\Pi\lt(m(x+h),m(x)\rt) \, dx,
	\end{align*} 
	where $C_0>0$ is an absolute constant.
\end{lem}
\begin{proof}
	We know from Lemma \ref{l:cENT} that if $m$ is $BV$ then 
	\begin{align*}
		\left(\bigvee_{\|\lambda_{\Phi}'\|_{\infty}\leq 1} \lt|\na \cdot \Phi(m)\rt|\right)
		\lt(\Omega'\rt)=\int_{J_m\cap \Omega'} c^{\mathrm{ENT}}\lt(m^{+},m^{-}\rt) \,d\mathcal{H}^1,
	\end{align*} 
	where $c^{\mathrm{ENT}}$ is defined by \eqref{eq:cENT}. Thanks to the inequality $\mathrm{c}^{\mathrm{ENT}}\leq\Pi$ provided by Lemma~\ref{l:cENTleqLambda}, we deduce
	\begin{align*}
		\left(\bigvee_{\|\lambda_{\Phi}'\|_{\infty}\leq 1} \lt|\na \cdot \Phi(m)\rt|\right)
		\lt(\Omega'\rt)\leq
		\int_{J_m\cap \Omega'} \Pi\lt(m^{+},m^{-}\rt) \,d\mathcal{H}^1.
	\end{align*}
	Hence the proof of Lemma~\ref{misbv} follows from the inequality
	\begin{align}
		\label{eq:intJmPi}
		\int_{J_m\cap \Omega'} \Pi\lt(m^{+},m^{-}\rt) \,d\mathcal{H}^1 \leq C_0 \, 
		\limsup_{\lt|h\rt|\to 0} \lt|h\rt|^{-1}\int_{\Omega'} 
		\Pi\lt(m(x+h),m(x)\rt) \, dx.
	\end{align}
	This inequality is valid for any Lipschitz function $\Pi$ and $BV$ map $m$, as a consequence of the rectifiability  of $J_m$ and the trace properties of $m$ (see e.g. \cite{ambrosio}), via  a Besicovitch covering argument which we detail next.

	Let $\delta\in (0,1)$. There exists $\ep_0>0$ and a subset $\wt{J}_m\subset J_m$ with $\mathcal{H}^1\lt(J_m\cap \Omega' \setminus \wt{J}_m\rt)<\delta$ and $\widetilde J_m +B_{\e_0}(0)\subset \Omega'$, such that for any $x_0\in \wt{J}_m$ and $0<r<\e_0$,
	\begin{align}
		\label{eqmisbv11}
		\Xint{-}_{B_r(x_0)\cap J_m} \lt|m^{\pm}\lt(x\rt)-m^{\pm}(x_0)\rt| d\mathcal{H}^1 (x)&<\delta,
		\nonumber\\
		\lt|\mathcal{H}^1\lt(B_r(x_0)\cap J_m\rt)-2r \rt|&<\delta \, r,\nn\\
		\text{ and } \Xint{-}_{B^{\pm}_r(x_0)} \lt|m\lt(x\rt)-m^{\pm}(x_0)\rt| dx &<\delta\quad\text{ for all }0<r<\ep_0,
	\end{align}
	where $B^{\pm}_r(x)$ denote the two half balls obtained by splitting $B_r(x)$ along the tangent line to $J_m$ at $x$.  Let $\ep\in (0,\ep_0)$. By Besicovitch's covering theorem \cite[Theorem~2.18]{ambrosio} there 
	exists an absolute constant $Q\in \mathbb{N}$ and families $\BI_1, \BI_2,\dots, \BI_Q$ of pairwise disjoint balls in the set $\lt\{B_{\ep}(x):x\in \wt{J}_m\rt\}$ such that 
	\begin{equation*}
		\wt{J}_m\subset \bigcup_{k=1}^Q \bigcup_{B\in \BI_k} B.
	\end{equation*}
	In particular for some $k_0\in \lt\{1,2,\dots, Q\rt\}$ we have 
	\begin{equation}\label{eq:intPileqQsumPij}
		\int_{\wt{J}_m} \Pi\lt(m^{+},m^{-}\rt) d\HI^1 \leq Q \sum_{B\in \BI_{k_0}} \int_{\wt{J}_m\cap B} \Pi\lt(m^{+},m^{-}\rt) d\HI^1.
	\end{equation}
	We have $\BI_{k_0}=\lbrace B_\e(x_j)\rbrace_{j=1,\ldots,p}$ for some $x_1,\ldots ,x_p\in\widetilde J_m$. 
	Using that $\Pi$ is
	Lipschitz thanks to its definition \eqref{eq:Lambda}, and the properties \eqref{eqmisbv11} of $\widetilde J_m$, we find
	\begin{align*}
		\int_{\wt{J}_m\cap B_{\ep}(x_j)} \Pi\lt(m^{+},m^{-}\rt) d\HI^1
		& \leq 
		2\ep   \Pi\lt(m^{+}(x_j),m^{-}(x_j)\rt) +  2L\,\delta \,\mathcal{H}^1\lt(J_m\cap B_{\ep}(x_j)\rt)\\
		&\leq 
		2\e
		\Xint{-}_{B^{+ }_{\ep}(x_j)}\Xint{-}_{B^{- }_{\ep}(x_j)}
		\Pi\lt(m(x), m\lt(\ti{x}\rt)\rt) \,d\ti{x}\, dx
		+5 L\delta\e,
	\end{align*}
	for some $L>0$ depending only on $\Pi$. Summing over $j=1,\ldots,p$ and taking \eqref{eq:intPileqQsumPij} into account, we deduce
	\begin{align*}
		\int_{\wt{J}_m} \Pi\lt(m^{+},m^{-}\rt) d\HI^1
		&\leq 2 Q\e 
		\sum_{j=1}^p \Xint{-}_{B^{+ }_{\ep}(x_j)}\Xint{-}_{B^{- }_{\ep}(x_j)}
		\Pi\lt(m(x), m\lt(\ti{x}\rt)\rt) \,d\ti{x}\, dx +5 L\delta\, p\e.
	\end{align*}
	Noting from the properties \eqref{eqmisbv11} of $\widetilde J_m$ that 
	\begin{align*}
		\mathcal{H}^1\lt(J_m\cap \Omega'\rt)\geq \sum_{k=1}^p \mathcal{H}^1\lt(B_{\ep}(x_k)\cap J_m\rt)\geq p \ep,
	\end{align*}
	this implies
	\begin{align*}
		\int_{\wt{J}_m} \Pi\lt(m^{+},m^{-}\rt) d\HI^1
		&\leq 2 Q\e 
		\sum_{j=1}^p \Xint{-}_{B^{+ }_{\ep}(x_j)}\Xint{-}_{B^{- }_{\ep}(x_j)}
		\Pi\lt(m(x), m\lt(\ti{x}\rt)\rt) \,d\ti{x}\, dx\\
		&\quad +5 L\delta\, \mathcal{H}^1\lt(J_m\cap \Omega'\rt).
	\end{align*}
	Moreover we have
	\begin{align*}
		&
		\e 
		\sum_{j=1}^p \Xint{-}_{B^{+ }_{\ep}(x_j)}\Xint{-}_{B^{- }_{\ep}(x_j)}
		\Pi\lt(m(x), m\lt(\ti{x}\rt)\rt) \,d\ti{x}\, dx \\
		&\leq 
		\frac{16}{\pi\e}
		\sum_{j=1}^p \int_{B^{+ }_{\ep}(x_j)}\left(\Xint{-}_{B_{2\ep}(0)}
		\Pi\lt(m(x), m\lt(x+h\rt)\rt)\,dh\right)\, dx \\
		&\leq 
		\frac{16}{\pi\e} \int_{\Omega'}\Xint{-}_{B_{2\ep}(0)}
		\Pi\lt(m(x), m\lt(x+h\rt)\rt)\,dh \, dx\\
		&\leq \frac{32}{\pi} 
		\sup_{|h|<2\e}
		\frac{1}{|h|}
		\int_{\Omega'}
		\Pi\lt(m(x), m\lt(x+h\rt)\rt) \, dx,
	\end{align*}
	provided $\ep<1/2\dist(\Omega',\partial\Omega)$, so plugging this into the previous inequality we deduce
	\begin{align*}
		\int_{\wt{J}_m} \Pi\lt(m^{+},m^{-}\rt) d\HI^1
		&
		\leq \frac{64 Q}{\pi}
		\sup_{|h|<2\e}
		\frac{1}{|h|}
		\int_{\Omega'}
		\Pi\lt(m(x), m\lt(x+h\rt)\rt) \, dx 
		+ 5 L \delta \, \mathcal{H}^1\lt(J_m\rt).
	\end{align*}
	Taking the limits $\e\to0$ and then $\delta\to 0$, we obtain \eqref{eq:intJmPi}, which concludes the proof of Lemma~\ref{misbv}.
\end{proof}

%
%
%
%

\bibliographystyle{acm}
\bibliography{aviles_giga}

\end{document}